\documentclass[12pt]{amsart}

\title{Monoidal model structures over infinite groups}
\author{Ioannis Emmanouil and Olympia Talelli}
\thanks{Research supported by the Hellenic Foundation
for Research and Innovation (H.F.R.I.) under the "3rd
Call for H.F.R.I.\ Research Projects to Support Faculty
Members and Researchers", project number 24921}

\newtheorem{Lemma}{Lemma}[section]
\newtheorem{Proposition}[Lemma]{Proposition}
\newtheorem{Theorem}[Lemma]{Theorem}
\newtheorem{Corollary}[Lemma]{Corollary}

\begin{document}

\begin{abstract}
The stable category of modules over the algebra of
a finite group with coefficients in a field is a
compactly generated tensor triangulated category, that
has been studied extensively in representation theory.
In this paper, we provide a plethora of infinite groups
$G$, for which the category of $kG$-modules (where $k$
is a commutative coherent ring of finite global dimension)
admits a monoidal model structure, in the sense of Hovey,
whose associated homotopy category is a compactly generated
tensor triangulated category. To that end, we use a
technique recently introduced by the authors, which is
based on Kropholler's operation ${\scriptstyle{{\bf LH}}}$
and the second author's operation $\Phi$.
\end{abstract}

\maketitle
\tableofcontents

\addtocounter{section}{-1}
\section{Introduction}

\noindent
A model structure in a category provides the framework for a
systematic study of the homotopical properties of objects and
morphisms therein, by analogy to the standard topological setting.
In the case of abelian categories, one looks for model structures
that are compatible with the additional data. A method for
constructing such model structures $\mathcal{M}$ in an abelian
category $\mathcal{A}$ is provided by the fundamental work of
Hovey \cite{Ho2} and Gillespie \cite{Gil1}, through the so-called
Hovey triples $(\mathcal{C},\mathcal{W},\mathcal{F})$. These
are particular classes of objects in $\mathcal{A}$ that induce
two complete cotorsion pairs (see \cite{EJ2} for this notion)
and fully determine the associated model structure; in particular,
$\mathcal{C}$ is the class of cofibrant objects and $\mathcal{F}$
is the class of fibrant objects. If both cotorsion pairs are
hereditary, then the intersection $\mathcal{C} \cap \mathcal{F}$,
with the exact structure inherited by $\mathcal{A}$, is Frobenius
and the stable category $\underline{\mathcal{C} \cap \mathcal{F}}$,
with the triangulated structure defined therein by Happel \cite{Ha},
is triangulated equivalent to the homotopy category
$\mbox{Ho}(\mathcal{M})$ of $\mathcal{M}$. In particular, if
$\mathcal{A}$ is a category with enough projectives, then any
complete hereditary projective cotorsion pair in $\mathcal{A}$
induces a Hovey triple and hence an associated model structure
$\mathcal{M}$ on $\mathcal{A}$, where all objects are fibrant.

Gorenstein homological algebra is the relative homological
theory which is based upon the classes of Gorenstein
projective, Gorenstein flat and Gorenstein injective modules,
introduced in \cite{EJ1} and \cite{EJT}. Having found interesting
applications in algebraic geometry, representation theory and
group cohomology, the theory has developed rapidly during the past
few decades. The class ${\tt GProj}(R)$ of Gorenstein projective
modules over a ring $R$ defines a hereditary projective cotorsion
pair $\left( {\tt GProj}(R),{\tt GProj}(R)^{\perp} \right)$ in the
module category $R$-Mod; see \cite{CIS}. The completeness of this
cotorsion pair though is only known in some cases. If $R$ is a left
Gorenstein ring, as defined by Beligiannis in \cite{Bel1}, then the
above cotorsion pair is complete; this has been proved by Holm in
\cite[Theorem 2.10]{Hol}. Bravo, Gillespie and Hovey have defined, 
over any ring $R$, a special type of Gorenstein projective modules, 
the so-called Gorenstein AC-projective modules, that provide a 
complete hereditary projective cotorsion pair and hence induce 
the Gorenstein AC-projective model structure; see 
\cite[Theorem 8.5]{BGH}. Another subclass of ${\tt GProj}(R)$,
which actually contains all Gorenstein AC-projective modules, was
introduced by \v{S}aroch and \v{S}t$\!$'$\!$ov\'{i}\v{c}ek 
\cite{SS}; this is the class of projectively coresolved Gorenstein
flat modules. The latter class defines a complete hereditary
projective cotorsion pair as well and hence induces the projectively
coresolved Gorenstein flat model structure; see the discussion
following \cite[Corollary 4.12]{SS}.

We are interested in the special case where $R$ is the algebra 
of a group $G$ with coefficients in a commutative ring $k$. Then, 
the Hopf algebra structure on $kG$ endows $kG$-Mod with the 
structure of a symmetric closed monoidal category. If $k$ is a 
field and $G$ a finite group, then the category of $kG$-modules
is Frobenius (with projective-injective objects the projective
modules) and the associated stable category $\underline{kG\mbox{-Mod}}$
is a compactly generated tensor triangulated category. It is a central
object of study in the (modular) representation theory of finite groups
and the corresponding literature is vast. Here, as an indication only, 
we mention the works \cite{BCR}, \cite{BIK}, \cite{BIKP1} and \cite{BIKP2}. 
The first attempt to construct model structures in the category of 
modules over infinite groups (or, at least, in suitable subcategories 
of the module category) is due to Benson \cite{Ben}. Kendall \cite{Ke} 
has recently shown that the Gorenstein projective model structure has 
a compactly generated tensor triangulated homotopy category over a big 
class of locally hierarchically defined groups, that contains the class 
${\scriptstyle{{\bf LH}}}\mathfrak{F}$ of groups defined by Kropholler 
in \cite{Kro}. In this paper, we examine the relevant homological conditions 
on $kG$-modules and study their behaviour with respect to certain 
operations on group classes. As an application of our analysis, we 
enlarge further the class of groups over which the Gorenstein 
projective model structure has a compactly generated tensor 
triangulated homotopy category.

The motivation behind the introduction of
${\scriptstyle{{\bf LH}}}\mathfrak{F}$-groups was the study
of certain homological finiteness conditions. Groups in
${\scriptstyle{{\bf LH}}}\mathfrak{F}$ are locally constructed
hierarchically, starting with the class $\mathfrak{F}$ of finite
groups; the resulting class of groups is very big and it does
require some effort to find groups which are not contained therein.
More generally, one may apply Kropholler's construction with base
any group class $\mathfrak{C}$ and obtain a group class
${\scriptstyle{{\bf LH}}}\mathfrak{C}$, which is (in principle,
at least) much bigger that $\mathfrak{C}$. Another example of an
operation on group classes is provided by Talelli's introduction
of groups of type $\Phi$, as a means to construct groups that admit
a finite dimensional model for the classifying space for proper
actions; see \cite{T}. The definition of groups of type $\Phi$
is based on the effort to describe modules of finite projective
dimension over the group, in terms of the finiteness of the
projective dimension of their restriction to finite subgroups.
More generally, starting with any group class $\mathfrak{C}$,
Talelli's definition provides a group class $\Phi\mathfrak{C}$;
see also \cite{MS}. Replacing the projective dimension with the
flat dimension in the discussion above, we obtain the group class
operation $\Phi_{flat}$.

Often, an iteration of Kropholler's operation ${\scriptstyle{{\bf LH}}}$
and Talelli's operations $\Phi$ and $\Phi_{flat}$ enlarges the class
of available examples of groups that satisfy appropriate homological
conditions. This is a technique introduced in \cite{EmmT3} and also
used in \cite{EmmR}. In this paper, we examine the following questions
on a group $G$, having fixed the commutative coefficient ring $k$:
(a) Assuming that the Gorenstein projective model structure is
defined in the full category of $kG$-modules, is it compatible
with the closed symmetric monoidal structure therein? (b) Is the
stable category of Gorenstein projective $kG$-modules compactly
generated? We prove in Sections 2 and 4 that the operations
${\scriptstyle{{\bf LH}}}$, $\Phi$ and $\Phi_{flat}$ behave
smoothly with respect to the homological conditions that answer
affirmatively questions (a) and (b) above. As an application,
let $\overline{\mathfrak{F}}$ be the class consisting of those
groups that are obtained from the class $\mathfrak{F}$ of finite
groups by a transfinite iteration process, which is based on the
group class operation
$\mathfrak{C} \mapsto {\scriptstyle{{\bf LH}}}\mathfrak{C} \cup
 \Phi\mathfrak{C} \cup \Phi_{flat}\mathfrak{C}$.
The following result is proved in Section 6; cf.\ Corollary 6.3.

\medskip

\noindent
{\bf Theorem.}
{\em Assume that the commutative ring $k$ is coherent and has
finite global dimension. Let $G$ be an
$\overline{\mathfrak{F}}$-group. Then, the Gorenstein projective
model structure in $kG$-Mod is defined and its homotopy category
is a compactly generated tensor triangulated category.}

\medskip

\noindent
If $k$ is a commutative coherent and Gorenstein regular ring
and $G$ is an $\overline{\mathfrak{F}}$-group, then there is
a generating set of compact objects for the stable category
of Gorenstein projective $kG$-modules, consisting of (Gorenstein
projective) $kG$-modules of type FP$_{\infty}$. This has the
following consequence: Over a pair $(k,G)$ as above, any 
finitely generated Gorenstein projective $kG$-module is
necessarily of type FP$_{\infty}$; cf.\ Corollary 5.2.

Here is a brief description of the contents of the paper:
Following Section 1, where certain basic preliminary notions
are recorded, in Section 2, we examine the behaviour of the
class of Gorenstein projective modules over a group algebra,
regarding tensor products. We examine properties implying
that the Gorenstein projective model structure is compatible
with the monoidal structure of the module category (i.e.\
that it is a monoidal model structure, in the sense of Hovey).
In Section 3, we detail certain more or less known properties
of the stable category of Gorenstein projective modules.
Then, in Section 4, we study the behaviour of the class
of groups over which the stable category of Gorenstein
projective modules is compactly generated, with respect
to the operations $\scriptstyle{\bf LH}$, $\Phi$ and
$\Phi_{flat}$. As an application, in Section 5, we
consider the finitely generated Gorenstein projective
modules over a group algebra and show that, over many
groups, these modules are necessarily of type
FP$_{\infty}$. In Section 6, we combine the results
of Sections 2 and 4 to prove the Theorem stated above.
Even though the class $\overline{\mathfrak{F}}$ is very
big, it is not exhaustive; in the Appendix, we exhibit
uncountably many finitely generated groups that are not
contained therein.

\medskip

\noindent
{\em Notations and terminology.}
All rings are assumed to be unital and associative. If $R$ is
a ring, then (unless otherwise specified) all $R$-modules are
assumed to be left $R$-modules. We denote by ${\tt Proj}(R)$ and
${\tt Flat}(R)$ the classes of projective and flat $R$-modules,
respectively.

\section{Preliminaries}

\noindent
In this preliminary section, we record certain basic notions
that are used in the sequel. These notions concern the relation
between cotorsion pairs and exact model structures, the hierarchical
structure of groups in the closure of group classes under certain
operations and the basic properties of Gorenstein modules over group
algebras.

\vspace{0.1in}

\noindent
{\sc I.\ Cotorsion pairs and model structures.}
If $R$ is a ring, then the $\mbox{Ext}^1$-pairing induces a
certain orthogonality relation between modules: If {\tt S} is
any class of modules, then the left (resp.\ the right) orthogonal
$^{\perp}{\tt S}$ (resp.\ ${\tt S}^{\perp}$) of ${\tt S}$ is
the class consisting of those modules $M$, for which
$\mbox{Ext}^1_R(M,S) = 0$ (resp.\ $\mbox{Ext}^1_R(S,M) = 0$)
for all $S \in {\tt S}$. If ${\tt U},{\tt V}$ are two module
classes, then we say that $({\tt U},{\tt V})$ is a cotorsion
pair (cf.\ \cite{EJ2}) if ${\tt U} = \! \, ^{\perp} {\tt V}$
and ${\tt U}^{\perp} = {\tt V}$. The cotorsion pair is hereditary
if $\mbox{Ext}^i_R(U,V)=0$ for all $i>0$ and all $U \in {\tt U}$
and $V \in {\tt V}$; it is complete if for any module $M$ there
exist short exact sequences
\[ 0 \longrightarrow V \longrightarrow U \longrightarrow M
     \longrightarrow 0
   \;\;\; \mbox{and } \;\;\;
   0 \longrightarrow M \longrightarrow V' \longrightarrow U'
     \longrightarrow 0 , \]
where $U,U' \in {\tt U}$ and $V,V' \in {\tt V}$. We say that
$({\tt U},{\tt V})$ is cogenerated by a class of modules
${\tt S} \subseteq {\tt U}$ if ${\tt V} = {\tt S}^{\perp}$.
A cotorsion pair which is cogenerated by a {\em set} of modules
is necessarily complete; cf.\ \cite{ET}. Finally, if the
intersection ${\tt U} \cap {\tt V}$ coincides with the class of
projective modules, we say that the cotorsion pair is projective.

A Hovey triple in $R$-Mod is a triple
$({\mathcal C} , {\mathcal W} , {\mathcal F})$ of module classes,
for which the two pairs
$({\mathcal C} , {\mathcal W} \cap {\mathcal F})$ and
$({\mathcal C} \cap {\mathcal W} , {\mathcal F})$ are complete
cotorsion pairs and the class ${\mathcal W}$ is thick (i.e.\
$\mathcal{W}$ is closed under direct summands and satisfies
the 2-out-of-3 property for short exact sequences). The work
of Gillespie \cite{Gil1}, which is itself based on the work
of Hovey \cite{Ho2}, provides a bijection between Hovey triples
and exact model structures in the module category; cf.\
\cite[Theorem 3.3]{Gil1}. Here, a model structure $\mathcal{M}$
on a category consists in specifying three classes of morphisms
therein, the fibrations, the cofibrations and the weak
equivalences, that have certain properties reminiscent of the
usual properties in the topological setting. These properties
imply that the formal inversion of the weak equivalences produces a
category $\mbox{Ho}(\mathcal{M})$, the so-called homotopy category
of the model structure $\mathcal{M}$, that has well-behaved and
controlled Hom-sets. For an excellent and comprehensive account
of this notion, which was originally introduced by Quillen in
\cite{Q}, the reader is referred to Hovey's book \cite{Ho1}.
As shown in \cite[Proposition 5.2]{Gil1}, for an exact model
structure $\mathcal{M}$ in $R$-Mod whose associated Hovey triple
$({\mathcal C} , {\mathcal W} , {\mathcal F})$ induces two
{\em hereditary} complete cotorsion pairs, the class
${\mathcal C} \cap {\mathcal F}$ is a Frobenius category with
projective-injective objects the modules in
${\mathcal C} \cap {\mathcal W} \cap {\mathcal F}$. According
to Happel \cite{Ha}, the stable category
$\underline{{\mathcal C} \cap {\mathcal F}}$, i.e.\
${\mathcal C} \cap {\mathcal F}$ modulo its projective-injective
objects ${\mathcal C} \cap {\mathcal W} \cap {\mathcal F}$, is
then triangulated. Finally, using
\cite[Proposition 4.4 and Corollary 4.8]{Gil1}, we infer that
the homotopy category $\mbox{Ho}(\mathcal{M})$ of such an exact
model structure $\mathcal{M}$ is triangulated equivalent to
$\underline{{\mathcal C} \cap {\mathcal F}}$.

\vspace{0.1in}

\noindent
{\sc II.\ Group class operations.}
An operation $\mathcal{O}$ on group classes assigns a group class
$\mathcal{O}\mathfrak{C}$ to each group class $\mathfrak{C}$, so
that $\mathfrak{C} \subseteq \mathcal{O}\mathfrak{C}$ for any group
class $\mathfrak{C}$ and
$\mathcal{O}\mathfrak{C}_1 \subseteq \mathcal{O}\mathfrak{C}_2$
for any two group classes $\mathfrak{C}_1,\mathfrak{C}_2$ with
$\mathfrak{C}_1 \subseteq \mathfrak{C}_2$. A group class
$\mathfrak{C}$ is said to be closed under an operation
$\mathcal{O}$ ($\mathcal{O}$-closed, for short) if
$\mathcal{O}\mathfrak{C} = \mathfrak{C}$. The $\mathcal{O}$-closure
$\overline{\mathfrak{C}}$ of a group class $\mathfrak{C}$ is the
smallest $\mathcal{O}$-closed group class containing $\mathfrak{C}$.
In this paper, we are mostly interested in the operation
${\scriptstyle{{\bf LH}}}$ defined by Kropholler in \cite{Kro} and
the operation $\Phi$ defined by the second author in \cite{T}.

If $\mathfrak{C}$ is any group class, then we define for each
ordinal $\alpha$ the group class
${\scriptstyle{{\bf H}}}_{\alpha}\mathfrak{C}$, using transfinite
induction. We let
${\scriptstyle{{\bf H}}}_0\mathfrak{C} = \mathfrak{C}$. For an
ordinal $\alpha >0$, the class
${\scriptstyle{{\bf H}}}_{\alpha}\mathfrak{C}$ consists of those
groups $G$ which admit a cellular action on a finite dimensional
contractible CW-complex $X$, in such a way that each isotropy
subgroup of the action belongs to
${\scriptstyle{{\bf H}}}_{\beta}\mathfrak{C}$ for some ordinal
$\beta < \alpha$. If $G,X$ are as above, $k$ is a commutative
ring and $M$ is a $kG$-module, then the cellular chain complex
of $X$ induces an exact sequence of $kG$-modules
\[ 0 \longrightarrow M_d \longrightarrow \ldots
     \longrightarrow M_1 \longrightarrow M_0
     \longrightarrow M \longrightarrow 0 , \]
where $d = \dim X$ and the $M_i$'s are direct sums of $kG$-modules
of the form $\mbox{ind}_H^G\mbox{res}_H^GM$, for some
${\scriptstyle{{\bf H}}}_{\beta}\mathfrak{C}$-subgroup $H$
of $G$ with $\beta < \alpha$. We say that a group belongs to
${\scriptstyle{{\bf H}}}\mathfrak{C}$ if it belongs to
${\scriptstyle{{\bf H}}}_{\alpha}\mathfrak{C}$ for some $\alpha$.
Finally, the class ${\scriptstyle{{\bf LH}}}\mathfrak{C}$ consists
of those groups $G$, over which any finitely generated subgroup
$H \subseteq G$ is contained in a suitable
${\scriptstyle{{\bf H}}}\mathfrak{C}$-subgroup $K=K(H) \subseteq G$.

The concept of groups of type $\Phi$, that were introduced by
Talelli in \cite{T}, corresponds to another operation, applied
to the class $\mathfrak{F}$ of finite groups. In general, if
$k$ is a commutative ring and $\mathfrak{C}$ is a class of groups,
then we define $\Phi\mathfrak{C}$ as the class consisting of those
groups $G$, over which a $kG$-module $M$ has finite projective dimension
if and only if there is an integer $n=n(M)$, such that
$\mbox{pd}_{kH}\mbox{res}_H^GM \leq n$ for any $\mathfrak{C}$-subgroup
$H \subseteq G$. See also \cite[Definition 5.6]{Bi}. We shall also
consider the version of the operation $\Phi$ for the flat dimension
of modules and define for any group class $\mathfrak{C}$ the
class $\Phi_{flat} \mathfrak{C}$ accordingly.

As shown in \cite[$\S $1.IV]{EmmT3}, the group operations
${\scriptstyle{{\bf LH}}}$, $\Phi$ and $\Phi_{flat}$ have an
additional property, termed as {\em continuity} therein. We
say that a group class operation $\mathcal{O}$ is continuous
provided that the following folds: Whenever
$\mathfrak{C}_{\alpha}$ is a group class defined for each
ordinal number $\alpha$, so that
$\mathfrak{C}_{\alpha} \subseteq \mathfrak{C}_{\beta}$ for all
$\alpha \leq \beta$, and $\mathfrak{C}$ is the class consisting
of those groups which are contained in $\mathfrak{C}_{\alpha}$ for
some ordinal number $\alpha$, then the group class
$\mathcal{O}\mathfrak{C}$ is the class consisting of those groups
which are contained in $\mathcal{O}\mathfrak{C}_{\alpha}$ for some
ordinal number $\alpha$.

The groups in the closure $\overline{\mathfrak{C}}$ of a class
$\mathfrak{C}$ under a continuous operation $\mathcal{O}$ admit
a simple hierarchical description, that we shall now explicit.
To that end, we define the classes $\mathfrak{C}_{\alpha}$ for
any ordinal number $\alpha$, using transfinite induction: We let
$\mathfrak{C}_0 = \mathfrak{C}$,
$\mathfrak{C}_{\alpha} = \mathcal{O}\mathfrak{C}_{\beta}$ if
$\alpha = \beta +1$ is a successor ordinal and
$\mathfrak{C}_{\alpha} = \bigcup_{\beta < \alpha}\mathfrak{C}_{\beta}$
if $\alpha$ is a limit ordinal. Then, $\overline{\mathfrak{C}}$ is the
class consisting of those groups which are contained in
$\mathfrak{C}_{\alpha}$, for some ordinal number $\alpha$.

We finally note that, if
$\mathcal{O}_1, \mathcal{O}_2, \ldots , \mathcal{O}_n$ are
(continuous) operations on group classes, then we may consider
the  (continuous) operation
$\mathcal{O}  = \mathcal{O}_1 \cup \mathcal{O}_2 \cup \ldots
 \cup \mathcal{O}_n$,
which is defined by letting
$\mathcal{O}\mathfrak{C} = \mathcal{O}_1\mathfrak{C} \cup
 \mathcal{O}_2\mathfrak{C} \cup \ldots \cup
 \mathcal{O}_n\mathfrak{C}$
for any group class $\mathfrak{C}$. The $\mathcal{O}$-closure
of $\mathfrak{C}$ is the smallest group class containing
$\mathfrak{C}$, that is closed under all of the operations
$\mathcal{O}_1, \mathcal{O}_2, \ldots , \mathcal{O}_n$.

\vspace{0.1in}

\noindent
{\sc III.\ Gorenstein projective and Gorenstein flat modules.}
Let $R$ be a ring. Then, an acyclic complex of projective modules
is called totally acyclic if it remains acyclic after applying the
functor $\mbox{Hom}_R(\_\!\_,P)$ for any projective module $P$. We
say that a module is Gorenstein projective if it is a cokernel of
a totally acyclic complex of projective modules; let ${\tt GProj}(R)$
be the class of these modules. An acyclic complex of flat modules
is called totally acyclic if it remains acyclic after applying
the functor $I \otimes_R\_\!\_$ for any injective right module
$I$. A module is Gorenstein flat if it is a cokernel of a totally
acyclic complex of flat modules. These notions were introduced by
Enochs, Jenda and Torecillas in \cite{EJ1} and \cite{EJT}; see
also \cite{Hol}. A distinguished class of Gorenstein flat modules
has been introduced by \v{S}aroch and \v{S}t$\!$'$\!$ov\'{i}\v{c}ek
in \cite{SS}: A module is projectively coresolved Gorenstein flat
if it is a cokernel of an acyclic complex of projective modules,
which remains acyclic after applying the functor $I \otimes_R\_\!\_$
for any injective right module $I$. As shown in
\cite[Theorem 4.4]{SS}, all projectively coresolved Gorenstein
flat modules are Gorenstein projective.

The Gorenstein projective dimension of a module $M$ is the length
of a shortest resolution of $M$ by Gorenstein projective modules.
Analogously, we define the Gorenstein flat dimension of $M$. The
rings over which the Gorenstein projective dimension of any module
is finite have been studied by Bennis and Mahbou \cite{BM}. These
are precisely the left Gorenstein regular rings, in the sense of
Beligiannis \cite{Bel1}; of course, they include all rings of finite
left global dimension. The rings over which all modules have finite
Gorenstein flat dimension have been studied by Christensen, Estrada
and Thompson \cite{CET}; we call these rings weakly Gorensrtein
regular.\footnote{There is no need to specify the side on which the
ring acts, since weak Gorenstein regularity is a left-right symmetric
property; this is demonstrated in \cite{CET}.} All rings of finite
weak global dimension are weakly Gorenstein regular. We also note
that any (left or right) Gorenstein regular ring is weakly Gorenstein
regular; see \cite[Theorem 3.7]{WYSZ} or \cite[Theorem 5.13]{CELTW}.
In the special case where the ring $R$ is isomorphic with its opposite
ring (for example, if $R$ is the group algebra of a group with
coefficients in a commutative ring), the latter result admits a
simple proof; cf.\ \cite[Remark 5.4(ii)]{Emm}.

\vspace{0.1in}

\noindent
{\sc IV.\ Gorenstein modules over group algebras.}
We consider the special case where  $R=kG$ is the algebra of
a group $G$ with coefficients in a commutative ring $k$. The
reader is referred to Brown's book \cite{Bro} for the basic
facts regarding the diagonal action of $G$ on the tensor product
$M \otimes_kN$ and the Hom-group $\mbox{Hom}_k(M,N)$ of two
$kG$-modules $M,N$. In this way, the category $kG$-Mod is
endowed with a symmetric closed monoidal structure. We note
that the tensor product $M \otimes_kN$ is a projective
$kG$-module, if $M$ is a projective $kG$-module and $N$ is a
$k$-projective $kG$-module. The reader is also referred to
\cite{Bro} for the basic properties of the restriction and
induction functors that are associated with a subgroup
$H \subseteq G$.

It is easily seen that induction maps ${\tt GProj}(kH)$ into
${\tt GProj}(kG)$ for any subgroup $H \subseteq G$; see, for
example, \cite[Lemma 2.6(i)]{EmmT1}. There are certain
conjectural properties of the class of Gorenstein projective
modules (over a group algebra) though that are only known to
hold in particular cases. A comprehensive class of groups,
over which Gorenstein projective modules have a smooth behaviour,
may be defined by letting $\mathfrak{X}$ be the class consisting
of those groups $G$, for which any acyclic complex of projective
$kG$-modules remains acyclic after applying the functor
$I \otimes_{kG} \_\!\_$ for any injective $kG$-module $I$. Hence,
$G$ is an $\mathfrak{X}$-group if all cokernels of acyclic complexes
of projective $kG$-modules are necessarily projectively coresolved
Gorenstein flat. This class contains all groups $G$, for which $kG$
is weakly Gorenstein regular. As shown in \cite{EmmT3}, the class
$\mathfrak{X}$ is subgroup-closed, ${\scriptstyle{{\bf LH}}}$-closed,
$\Phi$-closed and $\Phi_{flat}$-closed. For subsequent use, we record
the following properties of $\mathfrak{X}$-groups:
\newline
(1) If $H$ is an $\mathfrak{X}$-subgroup of $G$, then
restriction maps ${\tt GProj}(kG)$ into ${\tt GProj}(kH)$.
\newline
(2) If $G$ is an $\mathfrak{X}$-group, then the class
${\tt GProj}(kG)$ is closed under tensoring with any
$k$-projective $kG$-module.
\newline
(3) If $G$ is an $\mathfrak{X}$-group, then
$\left( {\tt GProj}(kG),{\tt GProj}(kG)^{\perp} \right)$ is
a hereditary projective cotorsion pair, which is cogenerated
by a set; in particular, it is complete.
\newline
(4) If $G$ is an $\mathfrak{X}$-group, then the orthogonal
class ${\tt GProj}(kG)^{\perp}$ contains all $kG$-modules
of finite flat dimension.
\newline
(5) If $G$ is an $\mathfrak{X}$-group, then any Gorenstein
projective $kG$-module of finite flat dimension is necessarily
projective.
\newline
Property (1) follows since the restriction of any acyclic complex
of projective $kG$-modules is an acyclic complex of projective
$kH$-modules, whereas (2) is a consequence of the fact that the
class of acyclic complexes of projective $kG$-modules remains
invariant under tensoring with any $k$-projective $kG$-module.
Regarding properties (3) and (4), these follow by invoking
\cite[Theorems 4.4 and 4.9]{SS}, since all Gorenstein projective
$kG$-modules are projectively coresolved Gorenstein flat over
an $\mathfrak{X}$-group $G$. Finally, property (5) follows
from property (4).

\section{Gorenstein properties of tensor products of $kG$-modules}

\noindent
Consider a commutative ring $k$ and a group $G$. For subsequent 
use in the paper, in this section we examine whether the tensor 
product of two $kG$-modules with the diagonal action of $G$ is 
Gorenstein projective or right orthogonal to the class of Gorenstein
projectives.

Being mainly interested in groups that are contained in the
class $\mathfrak{X}$ of $\S $1.IV, we can follow a standard
argument and prove the next result.

\begin{Proposition}
Assume that $k$ has finite weak global dimension and let $G$ be an
$\mathfrak{X}$-group. Then, the tensor product of two Gorenstein
projective $kG$-modules is Gorenstein projective.
\end{Proposition}
\vspace{-0.05in}
\noindent
{\em Proof.}
Our assumption on $k$ implies that any Gorenstein projective $kG$-module
is $k$-projective. This follows as in \cite[Proposition 1.1]{EmmT1} from
the result by Benson and Goodearl \cite{BG} and Neeman \cite{N2} on the
contractibility of acyclic complexes of projective $k$-modules that have
flat kernels; see also \cite{CH}. The result is therefore an immediate
consequence of property (2) in $\S $1.IV. \hfill $\Box$

\medskip

\noindent
In Section 6, we shall be interested in knowing whether the tensor
product of two $kG$-modules satisfies another (perhaps less natural)
condition. To be more precise, we shall need to know whether the right
Ext$^1$-orthogonal class ${\tt GProj}(kG)^{\perp}$ remains invariant
under tensoring with Gorenstein projective $kG$-modules. Let
$\mathfrak{A}$ be the class consisting of those $\mathfrak{X}$-groups
$G$ that have the following property: If $M,N$ are $kG$-modules with
$M \in {\tt GProj}(kG)$ and $N \in {\tt GProj}(kG)^{\perp}$, then
$M \otimes_kN \in {\tt GProj}(kG)^{\perp}$.

\begin{Lemma}
The class $\mathfrak{A}$ is subgroup-closed.
\end{Lemma}
\vspace{-0.05in}
\noindent
{\em Proof.}
Let $G$ be an $\mathfrak{A}$-group and consider a subgroup
$H \subseteq G$. We have to show that $H \in \mathfrak{A}$.
First of all, we note that $H$ is an $\mathfrak{X}$-group,
since $G \in \mathfrak{X}$ and $\mathfrak{X}$ is subgroup-closed.
Let $M,N$ be two $kH$-modules with
$M \in {\tt GProj}(kH)$ and $N \in {\tt GProj}(kH)^{\perp}$.
Then, $\mbox{ind}_H^GM \in {\tt GProj}(kG)$ and
$\mbox{coind}_H^GN \in {\tt GProj}(kG)^{\perp}$; the latter
claim follows from the Eckmann-Shapiro isomorphism and property
(1) in $\S $1.IV. Since $G \in \mathfrak{A}$, we have
$\mbox{ind}_H^GM \otimes_k \mbox{coind}_H^GN \in
 {\tt GProj}(kG)^{\perp}$
and hence
\[ \mbox{res}_H^G\mbox{ind}_H^GM \otimes_k
   \mbox{res}_H^G\mbox{coind}_H^GN =
   \mbox{res}_H^G \! \left( \mbox{ind}_H^GM \otimes_k
   \mbox{coind}_H^GN \right) \! \in {\tt GProj}(kH)^{\perp} . \]
We note that $M$ is a direct summand of $\mbox{res}_H^G\mbox{ind}_H^GM$
and $N$ is a direct summand of $\mbox{res}_H^G\mbox{coind}_H^GN$,
so that $M \otimes_kN$ is a direct summand of
$\mbox{res}_H^G\mbox{ind}_H^GM \otimes_k
 \mbox{res}_H^G\mbox{coind}_H^GN$.
The closure of the class ${\tt GProj}(kH)^{\perp}$ under direct
summands therefore implies that
$M \otimes_kN \in {\tt GProj}(kH)^{\perp}$. Thus, $H$ is an
$\mathfrak{A}$-group, as needed. \hfill $\Box$

\begin{Theorem}
The class $\mathfrak{A}$ is closed under the operations
${\scriptstyle{{\bf LH}}}$, $\Phi$ and $\Phi_{flat}$.
\end{Theorem}
\vspace{-0.05in}
\noindent
{\em Proof.}
We have to show that
$\mathfrak{A} = {\scriptstyle{{\bf LH}}}\mathfrak{A} \cup
 \Phi\mathfrak{A} \cup \Phi_{flat}\mathfrak{A}$,
i.e.\ that
${\scriptstyle{{\bf LH}}}\mathfrak{A} \cup \Phi\mathfrak{A}
 \cup \Phi_{flat}\mathfrak{A} \subseteq \mathfrak{A}$.
To that end, let $G$ be a group contained in
${\scriptstyle{{\bf LH}}}\mathfrak{A} \cup \Phi\mathfrak{A}
 \cup \Phi_{flat}\mathfrak{A}$.
First of all, we note that
\[ G \in {\scriptstyle{{\bf LH}}}\mathfrak{A} \cup \Phi\mathfrak{A}
   \cup \Phi_{flat}\mathfrak{A} \subseteq
   {\scriptstyle{{\bf LH}}}\mathfrak{X} \cup \Phi\mathfrak{X}
   \cup \Phi_{flat}\mathfrak{X} = \mathfrak{X} . \]
We now consider two $kG$-modules $M,N$ with $M \in {\tt GProj}(kG)$
and $N \in {\tt GProj}(kG)^{\perp}$. If $H$ is an $\mathfrak{A}$-subgroup
of $G$, then $\mbox{res}_H^GM \in {\tt GProj}(kH)$ (this is precisely
property (1) in $\S $1.IV) and
$\mbox{res}_H^GN \in {\tt GProj}(kH)^{\perp}$. Therefore, we
conclude that
\[ \mbox{res}_H^G \! \left( M \otimes_k N \right) \! =
   \mbox{res}_H^GM \otimes_k \mbox{res}_H^GN \in
   {\tt GProj}(kH)^{\perp} . \]
Since this is the case for any $\mathfrak{A}$-subgroup $H \subseteq G$,
we may apply \cite[Corollary 3.2]{EmmR} to the subgroup-closed subclass
$\mathfrak{A} \subseteq \mathfrak{X}$ and conclude that
$M \otimes_kN \in {\tt GProj}(kG)^{\perp}$. We have therefore
shown that $G \in \mathfrak{A}$, as needed. \hfill $\Box$

\vspace{0.1in}

\noindent
Assume that $k$ has finite weak global dimension and consider a group
$G$ for which the group algebra $kG$ is Gorenstein regular. Then, $G$
is an $\mathfrak{X}$-group, as we have noted in $\S $1.IV. Kendall has
proved in \cite{Ke} that $G$ is actually an $\mathfrak{A}$-group. Indeed,
in this case, any Gorenstein projective $kG$-module is $k$-projective
and the class ${\tt GProj}(kG)^{\perp}$ coincides with the class of all
$kG$-modules of finite projective dimension. Let $\mathfrak{G}$ be the
class consisting of all groups for which the group algebra $kG$ is
Gorenstein regular. Kendall has also shown that
${\scriptstyle{{\bf LH}}}\mathfrak{G} \subseteq \mathfrak{A}$. We
can enlarge $\mathfrak{G}$ and consider the class $\mathfrak{W}$
consisting of all groups $G$ for which the group algebra $kG$ is
weakly Gorenstein regular. Finally, we let $\overline{\mathfrak{W}}$
be the closure of $\mathfrak{W}$ under the operation
${\scriptstyle{{\bf LH}}} \cup \Phi \cup \Phi_{flat}$ and note that
$\overline{\mathfrak{W}}$-groups admit a hierarchical description,
as explained in $\S $1.II.

\begin{Corollary}
If $k$ has finite weak global dimension, then
$\overline{\mathfrak{W}}$ is a subclass of $\mathfrak{A}$.
\end{Corollary}
\vspace{-0.05in}
\noindent
{\em Proof.}
In view of Theorem 2.3, it suffices to show that
$\mathfrak{W} \subseteq \mathfrak{A}$. First of all,
$\mathfrak{W}$ is contained in $\mathfrak{X}$, as we have
noted in $\S $1.IV. Let $G$ be a $\mathfrak{W}$-group, so
that the group algebra $kG$ is weakly Gorenstein regular.
Then, the right orthogonal ${\tt GProj}(kG)^{\perp}$
coincides with the class of all $kG$-modules of finite
flat dimension; cf.\ \cite[Theorem 2.11]{DLW} or
\cite[Lemma 3.2]{WZ}. On the other hand, our assumption
on $k$ implies as in the proof of Proposition 2.1 that
any Gorenstein projective $kG$-module is $k$-projective.
Since the class of $kG$-modules of finite flat dimension
remains invariant under tensoring with any $k$-projective
(in fact, with any $k$-flat) $kG$-module, it follows readily
that $G \in \mathfrak{A}$. \hfill $\Box$

\vspace{0.1in}

\noindent
Assume that $k$ is weakly Gorenstein regular. Then, for any
finite group $G$ the group algebra $kG$ is weakly Gorenstein
regular as well; see, for example, \cite[Corollary 4.6]{HLGZ}.
It follows that $G \in \mathfrak{W}$, so that $\mathfrak{W}$
contains the class $\mathfrak{F}$ of finite groups. Recall
from the Introduction that $\overline{\mathfrak{F}}$ is the
closure of $\mathfrak{F}$ under the operation
${\scriptstyle{{\bf LH}}} \cup \Phi \cup \Phi_{flat}$.

\begin{Corollary}
If $k$ has finite weak global dimension, then
$\overline{\mathfrak{F}}$ is a subclass of $\mathfrak{A}$.
\end{Corollary}
\vspace{-0.05in}
\noindent
{\em Proof.}
In view of Corollary 2.4, it suffices to show that finite
groups are contained in $\mathfrak{W}$. This follows from
the discussion above, since a ring of finite weak global
dimension is certainly weakly Gorenstein regular. \hfill $\Box$

\section{The stable category of Gorenstein projective modules}

\noindent
In this Section, we detail certain (more or less known)
properties of the stable category of Gorenstein projective
modules, that will be needed later in the paper. We recall
the definition of the triangulated structure therein, we
examine the subcategory of Gorenstein projective modules
of type FP$_{\infty}$ and elaborate on the behaviour of
the induction functor that is associated with a subgroup
of a group.

\vspace{0.1in}

\noindent
{\sc I.\ The basics.}
Let $R$ be a ring. The category ${\tt GProj}(R)$ of Gorenstein
projective modules, with the exact structure inherited from the
full module category, is Frobenius with projective-injective
objects given by the projective modules; see, for example,
\cite[Proposition 2.2]{DEH}. We may therefore consider the
stable category $\underline{\tt GProj}(R)$ with the triangulated
structure defined by Happel \cite{Ha}. The general theory of
triangulated categories is presented in Neeman's book \cite{N1}.
Let $M,N$ be Gorenstein projective modules and
$f : M \longrightarrow N$ be an $R$-linear map. If
\begin{equation}
 0 \longrightarrow M \stackrel{\iota}{\longrightarrow} P
   \longrightarrow \Sigma M \longrightarrow 0
\end{equation}
is a short exact sequence of $R$-modules, where $P$ is projective
and $\Sigma M$ is Gorenstein projective, then the pushout $Z$
of $f$ and $\iota$ fits into a commutative diagram with exact rows
\[ \begin{array}{ccccccccc}
   0 & \longrightarrow & M & \stackrel{\iota}{\longrightarrow}
     & P & \longrightarrow & \Sigma M & \longrightarrow & 0 \\
     & & {\scriptstyle{f}} \downarrow & & \downarrow & & \parallel
     & & \\
   0 & \longrightarrow & N & \stackrel{\jmath}{\longrightarrow}
     & Z & \stackrel{p}{\longrightarrow} & \Sigma M
     & \longrightarrow & 0
   \end{array} \]
In this way, we obtain a so-called {\em standard} triangle
\begin{equation}
 M \stackrel{f}{\longrightarrow} N
   \stackrel{\jmath}{\longrightarrow} Z
   \stackrel{p}{\longrightarrow} \Sigma M
\end{equation}
in $\underline{\tt GProj}(R)$. The triangulated structure in
the stable category $\underline{\tt GProj}(R)$ is defined so
that the distinguished triangles therein are precisely those
triangles that are isomorphic (in $\underline{\tt GProj}(R)$)
to a standard triangle. Given a short exact sequence of Gorenstein
projective modules
\[ 0 \longrightarrow M \stackrel{f}{\longrightarrow} N
     \stackrel{g}{\longrightarrow} L \longrightarrow 0 , \]
it is easily seen that there is a triangle
\[ M \stackrel{f}{\longrightarrow} N
     \stackrel{g}{\longrightarrow} L
     \longrightarrow \Sigma M , \]
which is isomorphic in $\underline{\tt GProj}(R)$ to the standard
triangle (2); hence, the above triangle is distinguished in
$\underline{\tt GProj}(R)$. We conclude that $L$ is contained in
any triangulated subcategory of $\underline{\tt GProj}(R)$ that
contains both $M$ and $N$.

\begin{Lemma}
Let $\mathcal{T} \subseteq \underline{\tt GProj}(R)$ be a
triangulated subcategory and consider an exact sequence of
Gorenstein projective modules
\[ 0 \longrightarrow M_n \longrightarrow \ldots
     \longrightarrow M_1 \longrightarrow M_0
     \longrightarrow M \longrightarrow 0 , \]
for some $n \geq 0$. If $M_i \in \mathcal{T}$ for all
$i=0,1, \ldots ,n$, then $M \in \mathcal{T}$.
\end{Lemma}
\vspace{-0.05in}
\noindent
{\em Proof.}
The result is trivial if $n=0$ and follows from the discussion
above if $n=1$. In general, we use induction on $n$. Indeed, the
class ${\tt GProj}(R)$ is closed under kernels of epimorphisms
and hence the kernel of the map $M_0 \longrightarrow M$ is
Gorenstein projective as well. \hfill $\Box$

\vspace{0.1in}

\noindent
{\sc II.\ Gorenstein projective modules of type FP$_{\infty}$.}
An $R$-module $M$ is of type $\mbox{FP}_{\infty}$ if it
admits a resolution by finitely generated free $R$-modules. 
Bieri has proved that this condition on $M$ is equivalent to 
the requirement that the functors $\mbox{Ext}^n_R(M,\_\!\_)$ 
commute with filtered colimits for all $n \geq 0$; cf.\ 
\cite[Corollary 1.6]{Bie}. An immediate consequence of Bieri's 
criterion is that the class of $R$-modules of type FP$_{\infty}$ 
is thick; see also \cite[Theorem 1.7]{BP}. 

Let ${\tt GProj}^{fin}(R)$ denote the class of Gorenstein 
projective $R$-modules of type FP$_{\infty}$. These modules 
are precisely the cokernels of totally acyclic complexes of 
finitely generated free $R$-modules; see 
\cite[Theorem 4.2(vi)]{BDT}. In fact, any finitely generated 
Gorenstein projective $R$-module admits a coresolution by 
finitely generated free $R$-modules, all of whose cokernels are 
finitely generated Gorenstein projective modules. In particular,
if $M$ is a finitely generated Gorenstein projective module
there exists a short exact sequence as in (1), where $P$ is
finitely generated free and $\Sigma M$ is finitely generated
Gorenstein projective; it is clear that $\Sigma M$ is actually 
finitely presented. Moreover, if $M$ is Gorenstein projective 
of type FP$_{\infty}$, then $\Sigma M$ is also Gorenstein 
projective of type FP$_{\infty}$.

Let $\underline{\tt GProj}^{fin}(R)$ be the full subcategory
of the stable category $\underline{\tt GProj}(R)$ consisting
of those Gorenstein projective modules that are isomorphic (in
$\underline{\tt GProj}(R)$) to a module in ${\tt GProj}^{fin}(R)$,
i.e.\ to a Gorenstein projective module of type FP$_{\infty}$.

\begin{Lemma}
$\underline{\tt GProj}^{fin}(R)$ is a triangulated subcategory
of the stable category $\underline{\tt GProj}(R)$.
\end{Lemma}
\vspace{-0.05in}
\noindent
{\em Proof.}
Let $M$ be a Gorenstein projective $R$-module of type
FP$_{\infty}$. Then, as we noted above, $M$ is a cokernel
of an acyclic complex of finitely generated free
$R$-modules, all of whose cokernels are Gorenstein
projective of type FP$_{\infty}$. It follows that
$\Sigma^nM \in \underline{\tt GProj}^{fin}(R)$ for all
$n \in \mathbb{Z}$.

We now consider a distinguished triangle in
$\underline{\tt GProj}(R)$
\[ M \stackrel{f}{\longrightarrow} N
     \longrightarrow Z \longrightarrow \Sigma M , \]
where $M,N \in \underline{\tt GProj}^{fin}(R)$. We may
assume that this is the standard triangle associated
with a linear map $f : M \longrightarrow N$ between
two Gorenstein projective $R$-modules $M,N$ of type
FP$_{\infty}$. We may choose a short exact sequence
as in (1), where $P$ is finitely generated free and
$\Sigma M$ is Gorenstein projective of type
FP$_{\infty}$. The construction in the beginning of
$\S $3.I describes $Z$ as an extension of $\Sigma M$
by $N$. In view of the closure of the class of Gorenstein
projective modules of type FP$_{\infty}$ under extensions,
we conclude that $Z$ is a Gorenstein projective module
of type FP$_{\infty}$ as well. It follows that
$Z \in \underline{\tt GProj}^{fin}(R)$, as needed.
\hfill $\Box$

\vspace{0.1in}

\noindent
We are interested in the thick closure
${\tt K}(R) =
 \mbox{thick} \! \left[ \underline{\tt GProj}^{fin}(R) \right]$
of the triangulated subcategory $\underline{\tt GProj}^{fin}(R)$
in the stable category $\underline{\tt GProj}(R)$; ${\tt K}(R)$
is the full subcategory of $\underline{\tt GProj}(R)$, consisting
of those Gorenstein projective modules that are isomorphic (in
$\underline{\tt GProj}(R)$) to a direct summand of a Gorenstein
projective module of type FP$_{\infty}$. In other words, a
Gorenstein projective module $M$ is contained in ${\tt K}(R)$
if and only if there are Gorenstein projective modules $M'$ and
$N$, with $N$ of type FP$_{\infty}$, such that $M \oplus M'$ is
isomorphic with $N$ in $\underline{\tt GProj}(R)$. The
subcategory ${\tt K}(R) \subseteq \underline{\tt GProj}(R)$
is triangulated; in particular, if $M \in {\tt K}(R)$, then
$\Sigma M \in {\tt K}(R)$.

\begin{Lemma}
If $M$ is a finitely generated Gorenstein projective $R$-module
contained in ${\tt K}(R)$, then $M$ is of type FP$_{\infty}$.
\end{Lemma}
\vspace{-0.05in}
\noindent
{\em Proof.}
We begin with the special case where $M$ is finitely presented
and fix two Gorenstein projective $R$-modules $M'$ and $N$, with
$N$ of type FP$_{\infty}$, such that $M \oplus M'$ is isomorphic
with $N$ in $\underline{\tt GProj}(R)$. Then, there are linear
maps $f : M \oplus M' \longrightarrow N$ and
$g : N \longrightarrow M \oplus M'$, such that
$1_{M \oplus M'} - g \circ f :
 M \oplus M' \longrightarrow M \oplus M'$
factors through a projective module $P$. It follows that
$M \oplus M'$ is a direct summand of $N \oplus P$. Using
Bieri's criterion \cite[Corollary 1.6]{Bie}, we conclude
that the functor
$\mbox{Ext}^i_R(N,\_\!\_) = \mbox{Ext}^i_R(N \oplus P,\_\!\_)$
commutes with filtered colimits for all $i \geq 1$. Being a
direct summand of the latter functor, the functor
$\mbox{Ext}^i_R(M,\_\!\_)$ commutes with filtered colimits for
all $i \geq 1$ as well. Since $M$ is finitely presented, the
functor $\mbox{Hom}_R(M,\_\!\_)$ also commutes with filtered
colimits. Invoking Bieri's criterion once more, it follows
that $M$ is of type FP$_{\infty}$.

We now consider the general case and assume that $M$ is a
finitely generated Gorenstein projective module contained
in ${\tt K}(R)$. We fix a short exact sequence (1), where
$P$ is finitely generated free and $\Sigma M$ is finitely
presented Gorenstein projective. As we have noted above,
our assumption that $M$ is contained in ${\tt K}(R)$
implies that $\Sigma M$ is contained in ${\tt K}(R)$ as
well. Then, the first part of the proof implies that
$\Sigma M$ is of type FP$_{\infty}$. Since $P$ is clearly
of type FP$_{\infty}$, the closure of the class of modules
of type FP$_{\infty}$ under kernels of epimorphisms implies
that $M$ is also of type FP$_{\infty}$, as needed.
\hfill $\Box$

\vspace{0.1in}

\noindent
{\sc III.\ Induction from subgroups.}
The stable category $\underline{\tt GProj}(R)$ has coproducts.
Indeed, let $(M_i)_i$ be any family of Gorenstein projective
modules and consider the direct sum
$M = \bigoplus_iM_i \in {\tt GProj}(R)$ with the canonical
inclusions $\nu_i : M_i \longrightarrow M$. Then, for any
Gorenstein projective module $N$ the isomorphism
\[ \mbox{Hom}_{R}(M,N) \stackrel{\sim}{\longrightarrow}
   {\textstyle{\prod_i}} \mbox{Hom}_{R}(M_i,N) \]
maps any $f \in \mbox{Hom}_{R}(M,N)$ onto
$(f \circ \nu_i)_i \in \prod_i \mbox{Hom}_{R}(M_i,N)$. It is clear
that $f : M \longrightarrow N$ factors through a projective module
if and only if $f \circ \nu_i : M_i \longrightarrow N$ factors through
a projective module for all $i$. Therefore, by passage to the quotient,
we obtain an isomorphism
\[ \underline{\mbox{Hom}}_{R}(M,N) \stackrel{\sim}{\longrightarrow}
   {\textstyle{\prod_i}} \underline{\mbox{Hom}}_{R}(M_i,N) , \]
which maps any $[f] \in \underline{\mbox{Hom}}_{R}(M,N)$ onto
$(]f] \circ [\nu_i])_i = ([f \circ \nu_i])_i \in
 \prod_i \underline{\mbox{Hom}}_{R}(M_i,N)$.
Hence, $M$ is the coproduct of the family $(M_i)_i$ in the stable
category $\underline{\tt GProj}(R)$.

Since the stable category $\underline{\tt GProj}(R)$ has coproducts,
we may consider and look for compact objects therein. Recall that an
object $C$ in a triangulated category with coproducts $\mathcal{T}$
is compact if the functor
\[ \mbox{Hom}_{\mathcal{T}}(C,\_\!\_) : {\mathcal T}
   \longrightarrow \mbox{Ab} \]
preserves coproducts. We are interested in the special case where
$R=kG$ is the group algebra of a group $G$ with coefficients in a
commutative ring $k$ and fix such a pair $(k,G)$. A convenient
method for constructing compact objects in $\underline{\tt GProj}(kG)$
is by induction from suitable subgroups of $G$. To explain this, let
$H \subseteq G$ be a subgroup and assume that the restriction functor
$\mbox{res}_H^G$ maps ${\tt GProj}(kG)$ into ${\tt GProj}(kH)$. Since
that functor is exact and maps projective $kG$-modules into projective
$kH$-modules, it descends to a triangulated functor
\begin{equation}
 \mbox{res}_H^G : \underline{\tt GProj}(kG) \longrightarrow
 \underline{\tt GProj}(kH) .
\end{equation}
In view of the discussion above regarding coproducts in stable
categories of Gorenstein projective modules, it is clear that the
functor $\mbox{res}_H^G$ preserves coproducts. We also consider the
induction functor $\mbox{ind}_H^G$ and note that it maps Gorenstein
projective $kH$-modules into Gorenstein projective $kG$-modules.
Since that functor is also exact and maps projective $kH$-modules
into projective $kG$-modules, it too descends to a triangulated
functor
\begin{equation}
 \mbox{ind}_H^G : \underline{\tt GProj}(kH) \longrightarrow
 \underline{\tt GProj}(kG) .\footnote{This functor is defined
 under no additional assumptions on the subgroup $H \subseteq G$.}
\end{equation}
We now consider a module $M \in {\tt GProj}(kH)$ and a module
$N \in {\tt GProj}(kG)$. The isomorphism
\[ \mbox{Hom}_{kG} \! \left( \mbox{ind}_H^GM,N \right) \!
   \stackrel{\sim}{\longrightarrow}
   \mbox{Hom}_{kH} \! \left( M,\mbox{res}_H^GN \right) \]
maps any $f \in \mbox{Hom}_{kG} \! \left( \mbox{ind}_H^GM,N \right)$
onto the composition
$f \circ \eta \in \mbox{Hom}_{kH} \! \left( M,\mbox{res}_H^GN \right)$,
where $\eta : M \longrightarrow \mbox{res}_H^G\mbox{ind}_H^GM$ is the
canonical map (the unit of the adjunction between the module categories).
It is clear that $f : \mbox{ind}_H^GM \longrightarrow N$ factors through
a projective $kG$-module if and only if
$f \circ \eta : M \longrightarrow \mbox{res}_H^GN$ factors through
a projective $kH$-module. (This is a consequence of the fact that
projective modules are preserved by both restriction and induction.)
Therefore, by passage to the quotient, we obtain an isomorphism
\[ \underline{\mbox{Hom}}_{kG} \! \left( \mbox{ind}_H^GM,N \right) \!
   \stackrel{\sim}{\longrightarrow}
   \underline{\mbox{Hom}}_{kH} \! \left( M,\mbox{res}_H^GN \right) , \]
which maps any
$[f] \in \underline{\mbox{Hom}}_{kG} \! \left( \mbox{ind}_H^GM,N
 \right)$
onto
$[f] \circ [\eta] = [f \circ \eta] \in \underline{\mbox{Hom}}_{kH} \!
 \left( M,\mbox{res}_H^GN \right)$.
Hence, the pair $\left( \mbox{ind}_H^G , \mbox{res}_H^G \right)$
is an adjoint pair of functors between the stable categories of
Gorenstein projective modules.

\begin{Lemma}
Let $G$ be a group and consider a subgroup $H \subseteq G$. Assume
that the restriction of any Gorenstein projective $kG$-module is a
Gorenstein projective $kH$-module. Then, for any compact object
$C \in \underline{\tt GProj}(kH)$ the induced $kG$-module
$\mbox{ind}_H^GC$ is compact in $\underline{\tt GProj}(kG)$.
\end{Lemma}
\vspace{-0.05in}
\noindent
{\em Proof.}
In view of the discussion above, the triangulated functor
$\mbox{ind}_H^G$ in (4) is left adjoint to the coproduct
preserving functor $\mbox{res}_H^G$ in (3). Hence, the
result follows from general principles; see, for example,
\cite[Lemma 5.4.1(1)]{Kra1}. \hfill $\Box$

\vspace{0.1in}

\noindent
We now examine the behaviour of the (triangulated and
coproduct-preserving) functor $\mbox{ind}_H^G$ in (4)
with respect to the formation of localizing subcategories.
We recall that a subcategory $\mathcal{T}'$ of a triangulated
category with coproducts $\mathcal{T}$ is called localizing
if it is triangulated and closed under coproducts. If ${\tt S}$
is a set of objects in $\mathcal{T}$ and $\beta$ is an
uncountable cardinal, then $\mbox{Loc}^{\beta}({\tt S})$
is the smallest triangulated subcategory of $\mathcal{T}$
that contains ${\tt S}$ and is closed under the formation
of coproducts of fewer that $\beta$ objects. If
$\beta \leq \beta'$, then
$\mbox{Loc}^{\beta}({\tt S}) \subseteq
 \mbox{Loc}^{\beta'}({\tt S})$.
The union
$\mbox{Loc}({\tt S}) =
 \bigcup_{\beta} \mbox{Loc}^{\beta}({\tt S})$
is the smallest localizing subcategory of $\mathcal{T}$
that contains ${\tt S}$; it is the localizing subcategory
of $\mathcal{T}$ generated by ${\tt S}$.

\begin{Lemma}
Let $G$ be a group, $H \subseteq G$ a subgroup and ${\tt S}$
a set of Gorenstein projective $kH$-modules. If $M$ lies in
the localizing subcategory
$\mbox{Loc}({\tt S}) \subseteq \underline{\tt GProj}(kH)$
generated by ${\tt S}$, then
$\mbox{ind}_H^GM \in \mbox{Loc}({\tt S}') \subseteq
 \underline{\tt GProj}(kG)$,
where ${\tt S}' = \mbox{ind}_H^G{\tt S}$.
\end{Lemma}
\vspace{-0.05in}
\noindent
{\em Proof.}
Let $\beta$ be an uncountable regular cardinal with
$M \in \mbox{Loc}^{\beta}({\tt S})$. It suffices to
show that
\[ \mbox{ind}_H^G \! \left( \mbox{Loc}^{\beta}({\tt S})
   \right) \subseteq \mbox{Loc}^{\beta}({\tt S}') . \]
As shown in \cite[$\S $3.2]{N1}, the triangulated subcategory
$\mbox{Loc}^{\beta}({\tt S}) \subseteq \underline{\tt GProj}(kH)$
can be constructed using a transfinite process beginning
with ${\tt S}$, by successively closing with respect to
triangles and adding coproducts of fewer that $\beta$
objects. Since the functor $\mbox{ind}_H^G$ is triangulated
and preserves coproducts, it maps all objects obtained by
this transfinite process into any triangulated subgategory
of $\underline{\tt GProj}(kG)$ which contains ${\tt S}'$ and
is closed under coproducts of fewer that $\beta$ objects. In
particular, $\mbox{ind}_H^G$ maps all of these objects into
$\mbox{Loc}^{\beta}({\tt S}')$, as needed. \hfill $\Box$

\section{Compactly generated stable categories of Gorenstein projective modules}

\noindent
We fix a commutative ring $k$ and let $G$ be a group. Then,
the stable category of Gorenstein projective $kG$-modules
is a triangulated category with coproducts. In this section,
we examine whether that category is compactly generated.

In general, if $\mathcal{T}$ is a triangulated category with coproducts
and ${\tt C}$ is a set of compact objects therein, then the following
two conditions are equivalent:

(i) If $X$ is an object in $\mathcal{T}$ and
     $\mbox{Hom}_{\mathcal{T}}(\Sigma^nC,X)=0$ for all $n \in \mathbb{Z}$
     and $C \in {\tt C}$, then $X=0$.

(ii) $\mbox{Loc}({\tt C}) = \mathcal{T}$.
\newline
For a proof, the reader may consult \cite[Corollary 3.4.8]{Kra2}. We
say that $\mathcal{T}$ is compactly generated if there exists a set of
compact objects ${\tt C}$ satisfying these conditions; in that case,
${\tt C}$ is referred to as a set of compact generators for $\mathcal{T}$.

Let $\mathfrak{B}$ be the class of those $\mathfrak{X}$-groups
$G$, for which the stable category $\underline{\tt GProj}(kG)$
of Gorenstein projective $kG$-modules is compactly generated.

\begin{Proposition}
Let $\mathfrak{B}_0 \subseteq \mathfrak{B}$ be a subclass
and choose for any $H \in \mathfrak{B}_0$ a set ${\tt S}_H$
of compact generators of $\underline{\tt GProj}(kH)$. If $G$
is a group contained in
${\scriptstyle{{\bf LH}}}\mathfrak{B}_0 \cup
 \Phi\mathfrak{B}_0 \cup \Phi_{flat}\mathfrak{B}_0$,
then $\underline{\tt GProj}(kG)$ is compactly generated and
\[ {\tt T}_G = {\textstyle{\bigcup}} \left\{
   \mbox{ind}_H^G {\tt S}_H : H \subseteq G ,
   H \in \mathfrak{B}_0 \right\} \]
is a set of compact generators. In particular, $G \in \mathfrak{B}$.
\end{Proposition}
\vspace{-0.05in}
\noindent
{\em Proof.}
We note that any $\mathfrak{B}_0$-subgroup $H$ of $G$ is
contained in $\mathfrak{B} \subseteq \mathfrak{X}$ and hence
the restriction of any Gorenstein projective $kG$-module is
a Gorenstein projective $kH$-module; this is property (1) in
$\S $1.IV. We may therefore invoke Lemma 3.4 and conclude
that ${\tt T}_G$ is a set of compact objects in
$\underline{\tt GProj}(kG)$. We shall prove that ${\tt T}_G$
generates the stable category $\underline{\tt GProj}(kG)$, by
distinguishing three separate cases, according to whether:
(i) $G \in {\scriptstyle{{\bf LH}}}\mathfrak{B}_0$,
(ii) $G \in \Phi\mathfrak{B}_0$ and
(iii) $G \in \Phi_{flat}\mathfrak{B}_0$. The final clause in
the statement of the Proposition will then follow, since
$\mathfrak{B}_0 \subseteq \mathfrak{B} \subseteq \mathfrak{X}$
and hence
\[ G \in {\scriptstyle{{\bf LH}}}\mathfrak{B}_0 \cup
 \Phi\mathfrak{B}_0 \cup \Phi_{flat}\mathfrak{B}_0
 \subseteq
 {\scriptstyle{{\bf LH}}}\mathfrak{X} \cup
 \Phi\mathfrak{X} \cup \Phi_{flat}\mathfrak{X} = \mathfrak{X}. \]

(i) Assume that $G \in {\scriptstyle{{\bf H}}}\mathfrak{B}_0$
and proceed by induction on the ordinal $\alpha$, for which
$G \in {\scriptstyle{{\bf H}}}_{\alpha}\mathfrak{B}_0$. If
$\alpha =0$, then $G \in \mathfrak{B}_0$ and hence we may let
$H=G$ in the union defining ${\tt T}_G$, in order to conclude
that ${\tt S}_G \subseteq {\tt T}_G$. Since
$\underline{\tt GProj}(kG)$ is generated by ${\tt S}_G$, it
is a fortiori generated by ${\tt T}_G$. For the inductive step,
assume that $\alpha >0$ and the result is true for all
${\scriptstyle{{\bf H}}}_{\beta}\mathfrak{B}_0$-groups and all
$\beta < \alpha$. If $M$ is a Gorenstein projective $kG$-module,
then, as noted in $\S $1.II, there exists an exact sequence of
$kG$-modules
\[ 0 \longrightarrow M_d \longrightarrow \ldots
     \longrightarrow M_1 \longrightarrow M_0
     \longrightarrow M \longrightarrow 0 , \]
where $d$ is the dimension of the $G$-CW-complex witnessing
that $G \in {\scriptstyle{{\bf H}}}_{\alpha}\mathfrak{B}_0$
and each $M_i$ is a direct sum of $kG$-modules of the form
$\mbox{ind}_H^G\mbox{res}_H^GM$, with $H$ an
${\scriptstyle{{\bf H}}}_{\beta}\mathfrak{B}_0$-subgroup of
$G$ for some $\beta < \alpha$. Since such a group $H$ is
contained in
${\scriptstyle{{\bf H}}}_{\beta}\mathfrak{B}_0 \subseteq
 {\scriptstyle{{\bf H}}}\mathfrak{B}_0 \subseteq
 {\scriptstyle{{\bf H}}}\mathfrak{B} \subseteq
 {\scriptstyle{{\bf H}}}\mathfrak{X} = \mathfrak{X}$,
the $kH$-module $\mbox{res}_H^GM$ is Gorenstein projective;
this is property (1) in $\S $1.IV. We also note that
\begin{equation}
 \mbox{ind}_H^G{\tt T}_H = {\textstyle{\bigcup}}
 \left\{ \mbox{ind}_H^G\mbox{ind}_K^H{\tt S}_K :
 K \subseteq H , K \in \mathfrak{B}_0 \right\}
 \subseteq {\textstyle{\bigcup}} \left\{
 \mbox{ind}_K^G{\tt S}_K : K \subseteq G,
 K \in \mathfrak{B}_0 \right\} = {\tt T}_G
\end{equation}
for such an $H$. Our inductive assumption on $H$ implies
that $\mbox{res}_H^GM \in {\tt Loc}({\tt T_H})$ and hence
\[ \mbox{ind}_H^G\mbox{res}_H^GM \in
   \mbox{ind}_H^G \! \left( {\tt Loc}({\tt T_H}) \right)
   \subseteq
   {\tt Loc} \! \left( \mbox{ind}_H^G{\tt T_H} \right)
   \subseteq {\tt Loc}({\tt T_G}) ; \]
cf.\ Lemma 3.5. It follows that
$M_i \in {\tt Loc}({\tt T_G})$ for all $i=0,1, \ldots ,d$
and hence an application of Lemma 3.1 implies that
$M \in {\tt Loc}({\tt T_G})$. Since this is the case for
any Gorenstein projective $kG$-module $M$, we conclude
that ${\tt Loc}({\tt T_G}) = \underline{\tt GProj}(kG)$.
This completes the inductive step of the proof, taking care
of any ${\scriptstyle{{\bf H}}}\mathfrak{B}_0$-group.

We now assume that
$G \in {\scriptstyle{{\bf LH}}}\mathfrak{B}_0$ and consider
a Gorenstein projective $kG$-module $M$, which is non-zero
in the stable category $\underline{\tt GProj}(kG)$. Then,
the $kG$-module $M$ is not projective and hence it is not
flat either; see property (5) in $\S $1.IV. We now claim that
there exists a finitely generated subgroup $K \subseteq G$,
such that the $kK$-module $\mbox{res}_K^GM$ is not projective.
Indeed, if $\mbox{res}_K^GM$ was projective for any finitely
generated subgroup $K \subseteq G$, then we would also have
$\mbox{res}_K^GM \in {\tt Flat}(kK)$ and hence
$\mbox{ind}_K^G\mbox{res}_K^GM \in {\tt Flat}(kG)$. Since
$M$ is the filtered colimit of the family
$\left( \mbox{ind}_K^G\mbox{res}_K^GM \right) _K$, it would
then follow that $M$ is flat, which is absurd. Let $K$ be a
finitely generated subgroup of $G$, such that
$\mbox{res}_K^GM \notin {\tt Proj}(kK)$. Since
$G \in {\scriptstyle{{\bf LH}}}\mathfrak{B}_0$, the subgroup
$K$ is contained in an
${\scriptstyle{{\bf H}}}\mathfrak{B}_0$-subgroup $H \subseteq G$.
We note that $\mbox{res}_K^GM = \mbox{res}_K^H\mbox{res}_H^GM$
and hence the $kH$-module $\mbox{res}_H^GM$ is not projective.
On the other hand,
$H \in {\scriptstyle{{\bf H}}}\mathfrak{B}_0 \subseteq
 {\scriptstyle{{\bf H}}}\mathfrak{B} \subseteq
 {\scriptstyle{{\bf H}}}\mathfrak{X} = \mathfrak{X}$,
so that the $kH$-module $\mbox{res}_H^GM$ is Gorenstein
projective; this is property (1) in $\S $1.IV. It follows
that $\mbox{res}_H^GM \neq 0 \in \underline{\tt GProj}(kH)$.
As we have already shown that $\underline{\tt GProj}(kH)$ is
generated by ${\tt T}_H$, there must exist a suitable
$C \in {\tt T}_H$ and $n \in \mathbb{Z}$, such that the abelian
group
\[ \underline{\mbox{Hom}}_{kG} \! \left( \Sigma^n\mbox{ind}_H^GC,M
   \right) =
   \underline{\mbox{Hom}}_{kG} \! \left( \mbox{ind}_H^G \Sigma^nC,M
   \right) \! \stackrel{\sim}{\longrightarrow}
   \underline{\mbox{Hom}}_{kH} \! \left( \Sigma^nC,\mbox{res}_H^GM
   \right) \]
is non-trivial. Since
$\mbox{ind}_H^GC \in \mbox{ind}_H^G{\tt T}_H \subseteq {\tt T}_G$
(use the same reasoning as in (5) above) and $M$ was an arbitrary
non-zero object in the stable category $ \underline{\tt GProj}(kG)$,
it follows that $\underline{\tt GProj}(kG)$ is indeed generated by
${\tt T}_G$.

(ii) Let $G$ be a $\Phi\mathfrak{B}_0$-group and consider a
Gorenstein projective $kG$-module $M$, which is non-zero in
the stable category $\underline{\tt GProj}(kG)$. Then, $M$
is not projective. We claim that there exists a
$\mathfrak{B}_0$-subgroup $H \subseteq G$, such that the $kH$-module
$\mbox{res}_H^GM$ is not projective either. Indeed, if $\mbox{res}_H^GM$
was projective for any $\mathfrak{B}_0$-subgroup $H \subseteq G$, then
we would have $\mbox{pd}_{kG}M < \infty$. Since any Gorenstein projective
module of finite projective dimension is necessarily projective by
\cite[Proposition 2.27]{Hol}, it would follow that $M$ is projective,
which is again absurd. We may therefore find a $\mathfrak{B}_0$-subgroup
$H \subseteq G$, such that $\mbox{res}_H^GM \notin {\tt Proj}(kH)$.
On the other hand,
$H \in \mathfrak{B}_0 \subseteq \mathfrak{B} \subseteq \mathfrak{X}$
and hence the $kH$-module $\mbox{res}_H^GM$ is Gorenstein projective;
this is property (1) in $\S $1.IV. We conclude that $\mbox{res}_H^GM$
is a non-zero object in $\underline{\tt GProj}(kH)$. The latter category
being generated by ${\tt S}_H$, there exists a suitable
$C \in {\tt S}_H$ and $n \in \mathbb{Z}$, such that the abelian group
\[ \underline{\mbox{Hom}}_{kG} \! \left( \Sigma^n\mbox{ind}_H^GC,M
   \right) =
   \underline{\mbox{Hom}}_{kG} \! \left( \mbox{ind}_H^G \Sigma^nC,M
   \right) \! \stackrel{\sim}{\longrightarrow}
   \underline{\mbox{Hom}}_{kH} \! \left( \Sigma^nC,\mbox{res}_H^GM
   \right) \]
is non-trivial. Since
$\mbox{ind}_H^GC \in \mbox{ind}_H^G{\tt S}_H \subseteq {\tt T}_G$
and $M$ was an arbitrary non-zero object in the stable category
$\underline{\tt GProj}(kG)$, it follows that $\underline{\tt GProj}(kG)$
is indeed generated by ${\tt T}_G$.

(iii) The proof in the case where $G$ is a
$\Phi_{flat}\mathfrak{B}_0$-group is completely analogous to that
given in (ii) above; we simply use property (5) in $\S $1.IV instead
of \cite[Proposition 2.27]{Hol}. \hfill $\Box$

\begin{Theorem}
The class $\mathfrak{B}$ is closed under the operations
${\scriptstyle{{\bf LH}}}$, $\Phi$ and $\Phi_{flat}$.
\end{Theorem}
\vspace{-0.05in}
\noindent
{\em Proof.}
This is an immediate consequence of Proposition 4.1, by
letting $\mathfrak{B}_0 = \mathfrak{B}$ therein. \hfill $\Box$

\vspace{0.1in}

\noindent
Let $\mathfrak{B}_0 \subseteq \mathfrak{B}$ be any subclass
and consider its closure $\overline{\mathfrak{B}_0}$ under 
the operation
${\scriptstyle{{\bf LH}}} \cup \Phi \cup \Phi_{flat}$.
We note that $\overline{\mathfrak{B}_0}$-groups admit a
hierarchical description, as explained in $\S $1.II. In
view of Theorem 4.2, we have an inclusion
$\overline{\mathfrak{B}_0} \subseteq \mathfrak{B}$, i.e.\
for any $\overline{\mathfrak{B}_0}$-group $G$ the stable
category $\underline{\tt GProj}(kG)$ is compactly generated.
The following result generalizes Proposition 4.1.

\begin{Corollary}
Let $\mathfrak{B}_0 \subseteq \mathfrak{B}$ be a subclass
and choose for any $H \in \mathfrak{B}_0$ a set ${\tt S}_H$
of compact generators of the stable category
$\underline{\tt GProj}(kH)$. If $G$ is a group contained in
the closure $\overline{\mathfrak{B}_0}$ of $\mathfrak{B}_0$
under the operation
${\scriptstyle{{\bf LH}}} \cup \Phi \cup \Phi_{flat}$, then
$\underline{\tt GProj}(kG)$ is compactly generated and
\[ {\tt T}_G = {\textstyle{\bigcup}} \left\{
   \mbox{ind}_H^G {\tt S}_H : H \subseteq G ,
   H \in \mathfrak{B}_0 \right\} \]
is a set of compact generators.
\end{Corollary}
\vspace{-0.05in}
\noindent
{\em Proof.}
First of all, we note that property (1) in $\S $1.IV and Lemma
3.4 may be used as in the beginning of the proof of Proposition
4.1 to show that ${\tt T}_G$ is indeed a set of compact objects
in $\underline{\tt GProj}(kG)$. In order to describe the hierarchical
structure of $\overline{\mathfrak{B}_0}$-groups, we define the classes
$\mathfrak{B}_{0,\alpha}$ for any ordinal number $\alpha$, as follows:
\newline
(a) $\mathfrak{B}_{0,0} = \mathfrak{B}_0$,
\newline
(b) $\mathfrak{B}_{0,\alpha +1} =
     {\scriptstyle{{\bf LH}}}\mathfrak{B}_{0,\alpha}
     \cup \Phi\mathfrak{B}_{0,\alpha} \cup
     \Phi_{flat}\mathfrak{B}_{0,\alpha}$
    for any ordinal $\alpha$ and
\newline
(c) $\mathfrak{B}_{0,\alpha} =
     \bigcup_{\beta < \alpha}\mathfrak{B}_{0,\beta}$
    for any limit ordinal $\alpha$.
\newline
Then, $\overline{\mathfrak{B}_0}$ is the class consisting of
those groups $G$, for which $G \in \mathfrak{B}_{0,\alpha}$
for some ordinal $\alpha$. We may therefore prove the claim
in the statement of the Corollary using induction on the least
ordinal $\alpha$, for which  $G \in \mathfrak{B}_{0.\alpha}$.
As noted above, the class $\overline{\mathfrak{B}_0}$ is a
subclass of $\mathfrak{B}$; in particular, $\mathfrak{B}_{0,\alpha}$
is a subclass of $\mathfrak{B}$ for all $\alpha$.

In the cases where $\alpha =0$ or $\alpha =1$, the result follows
from Proposition 4.1. For the inductive step of the proof, we assume
that $G$ is a $\mathfrak{B}_{0,\alpha +1}$-group and the result
is known for all $\mathfrak{B}_{0,\alpha}$-groups. For any
$\mathfrak{B}_{0,\alpha}$-group $H$ our inductive assumption
implies that ${\tt T}_H$ is a set of compact generators for the
stable category $\underline{\tt GProj}(kH)$. Since
$\mathfrak{B}_{0,\alpha}$ is a subclass of $\mathfrak{B}$ and
\[ G \in \mathfrak{B}_{0,\alpha +1} =
   {\scriptstyle{{\bf LH}}}\mathfrak{B}_{0,\alpha}
   \cup \Phi\mathfrak{B}_{0,\alpha} \cup
   \Phi_{flat}\mathfrak{B}_{0,\alpha} , \]
we may apply Proposition 4.1 and conclude that a set of compact
generators for $\underline{\tt GProj}(kG)$ is provided by
\[ {\tt T}'_G = {\textstyle{\bigcup}} \left\{
   \mbox{ind}_H^G {\tt T}_H : H \subseteq G ,
   H \in \mathfrak{B}_{0,\alpha} \right\} . \]
As $\mbox{ind}_H^G {\tt T}_H \subseteq {\tt T}_G$ for any
$\mathfrak{B}_{0,\alpha}$-subgroup $H \subseteq G$ (see the
displayed equation (5) in the proof of Proposition 4.1), it
follows that ${\tt T}'_G \subseteq {\tt T}_G$. In particular,
${\tt T}_G$ is also a set of compact generators for the stable
category $\underline{\tt GProj}(kG)$, as needed. \hfill $\Box$

\vspace{0.1in}

\noindent
Following \cite{BGH}, we say that a $kG$-module $L$ is {\em level}
if the functor $\mbox{Tor}_1^{kG}(L, \_\!\_)$ vanishes on all
$kG$-modules of type $\mbox{FP}_{\infty}$. Flat modules are
clearly level; in fact, the class of level modules contains
the definable closure of the regular module. A $kG$-module is
called Gorenstein AC-projective if it is a cokernel of an acyclic
complex of projective $kG$-modules, which remains acyclic after
applying the functor $\mbox{Hom}_{kG}(\_\!\_ ,L)$ for any level
$kG$-module $L$. Here, the capital initials refer to the {\em
absolutely clean} modules; these are the modules in the right
Ext$^1$-orthogonal to the class of $kG$-modules of type
$\mbox{FP}_{\infty}$. All Gorenstein AC-projective modules are
Gorenstein projective; in fact, they are even projectively
coresolved Gorenstein flat (cf.\ \cite[Corollary 4.5]{SS}).
Following \cite{Ke}, let $\mathfrak{G}_{AC}$ be the class
consisting of those groups $G$, over which any $kG$-module
admits a resolution of finite length by Gorenstein 
AC-projective modules; in this case, we say that $kG$ has 
finite global Gorenstein AC-projective dimension. The class
$\mathfrak{G}_{AC}$ is contained in $\mathfrak{G}$ and hence
$\mathfrak{G}_{AC} \subseteq \mathfrak{X}$. We consider the
closure $\overline{\mathfrak{G}_{AC}}$ of $\mathfrak{G}_{AC}$
under the operation
${\scriptstyle{{\bf LH}}} \cup \Phi \cup \Phi_{flat}$ and
note that $\overline{\mathfrak{G}_{AC}}$-groups admit a
hierarchical description, as explained in $\S $1.II.

\begin{Corollary}
The class $\overline{\mathfrak{G}_{AC}}$ is contained in
$\mathfrak{B}$.
\end{Corollary}
\vspace{-0.05in}
\noindent
{\em Proof.}
In view of Theorem 4.2, it suffices to show the inclusion
$\mathfrak{G}_{AC} \subseteq \mathfrak{B}$. First of all,
$\mathfrak{G}_{AC}$ is contained in $\mathfrak{X}$, as we
have already noted above. Let $G$ be a $\mathfrak{G}_{AC}$-group.
Then, as shown in \cite[Lemma 2.21]{Ke}, all level $kG$-modules
have finite injective dimension. It follows that any Gorenstein
projective $kG$-module is Gorenstein AC-projective. It also follows
that $kG$ is an AC-Gorenstein ring, in the sense of Gillespie
\cite{Gil2}. Hence, \cite[Theorem 6.2]{Gil2} implies that the 
model structure in the category of $kG$-modules that is 
associated with the cotorsion pair
$\left( {\tt GProj}(kG),{\tt GProj}(kG)^{\perp} \right)$ has
a compactly generated homotopy category. Since the latter
homotopy category is triangulated equivalent to the stable
category $\underline{\tt GProj}(kG)$, we conclude that it
too is compactly generated. (For later use, we note that
\cite[Theorem 6.2]{Gil2} provides us with a generating set
of compact objects for the stable category
$\underline{\tt GProj}(kG)$, which consists of Gorenstein
projective $kG$-modules of type FP$_{\infty}$.) It follows
that $G \in \mathfrak{B}$. \hfill $\Box$

\vspace{0.1in}

\noindent
Recall that $\overline{\mathfrak{F}}$ is the closure of the
class $\mathfrak{F}$ of finite groups under the operation
${\scriptstyle{{\bf LH}}} \cup \Phi \cup \Phi_{flat}$.

\begin{Corollary}
If $k$ is coherent and Gorenstein regular, then
$\overline{\mathfrak{F}}$ is a subclass of $\mathfrak{B}$.
\end{Corollary}
\vspace{-0.05in}
\noindent
{\em Proof.}
In view of Corollary 4.4, it suffices to show that finite
groups are contained in $\mathfrak{G}_{AC}$. To that end,
let $G$ be a finite group. Since $k$ is Gorenstein regular,
the group algebra $kG$ is also Gorenstein regular; see, for
example, \cite[Corollary 3.4]{CR}. The coherence of $k$
implies that the group algebra $kG$ is coherent as
well.\footnote{Assume that $\mid \! G \! \mid \; = n$ and
let $I$ be a finitely generated left ideal of $kG$. We may
write $I \simeq F/K$, where $F$ is a finitely generated free
$kG$-module. Then, $\mbox{res}_1^GI$ is a finitely generated
submodule of $\mbox{res}_1^GkG = k^n$ and hence $\mbox{res}_1^GI$
is a finitely presented $k$-module. It follows that $K$ is
finitely generated as a $k$-module and a fortiori it is
finitely generated as a $kG$-module.} Hence, any finitely
presented $kG$-module is of type $\mbox{FP}_{\infty}$, so
that the level $kG$-modules are precisely the flat $kG$-modules.
Invoking \cite[Corollary 4.5]{SS}, we conclude that all
projectively coresolved Gorenstein flat $kG$-modules are
Gorenstein AC-projective. On the other hand, $k$ is weakly
Gorenstein regular and hence \cite[Proposition 3.4]{EmmT3}
implies that $G \in \mathfrak{X}$. It follows that any
Gorenstein projective $kG$-module is projectively coresolved
Gorenstein flat and hence Gorenstein AC-projective. Therefore,
the Gorenstein regularity of $kG$ implies that any $kG$-module
admits a finite length resolution by Gorenstein AC-projective
modules. Hence, $G \in \mathfrak{G}_{AC}$, as needed. \hfill $\Box$

\section{An application: Finitely generated Gorenstein projective modules}

\noindent
As an application of the results in Section 4, in this section
we examine the finitely generated Gorenstein projective modules
over the algebra of a group $G$ with coefficients in a commutative
ring $k$ and show that these modules are necessarily of type
FP$_{\infty}$, for many pairs $(k,G)$. 

We note that for any subgroup $H$ of a group $G$ the induction 
functor from the category of $kH$-modules to the category of 
$kG$-modules preserves the class of modules of type FP$_{\infty}$. 
In the special case of an 
${\scriptstyle{{\bf LH}}}\mathfrak{G}_{AC}$-group, the next result 
follows from \cite[Theorem 1.4]{Ke}.

\begin{Proposition}
Let $G$ be a a group contained in
$\overline{\mathfrak{G}_{AC}}$. Then, the following conditions
are equivalent for a $kG$-module $M$:

(i) $M$ is a cokernel of an acyclic complex of finitely
    generated free $kG$-modules,

(ii) $M$ is Gorenstein projective of type FP$_{\infty}$,

(iii) $M$ is Gorenstein projective and finitely presented and

(iv) $M$ is Gorenstein projective and finitely generated.
\end{Proposition}
\vspace{-0.05in}
\noindent
{\em Proof.}
We have already noted at the beginning of $\S $3.II that the
Gorenstein projective $kG$-modules of type FP$_{\infty}$ are
the cokernels of the totally acyclic complexes of finitely
generated free $kG$-modules. Since $\overline{\mathfrak{G}_{AC}}$
is a subclass of the class $\mathfrak{X}$ defined in $\S $1.IV,
any acyclic complex of projective $kG$-modules is necessarily
totally acyclic. This proves the equivalence (i)$\leftrightarrow$(ii).

It is clear that (ii)$\rightarrow$(iii)$\rightarrow$(iv)
and hence it only remains to show that (iv)$\rightarrow$(ii).
To that end, assume that $M$ is Gorenstein projective and
finitely generated. Then, $M$ is clearly a compact object of
the stable category $\underline{\tt GProj}(kG)$. The latter
category is compactly generated, in view of Corollary 4.4.
In fact, as we shall now explain, it admits a set of compact
generators that consists of Gorenstein projective $kG$-modules
of type FP$_{\infty}$. We have noted during the proof of Corollary
4.4 that for any ${\mathfrak{G}_{AC}}$-subgroup $H \subseteq G$
the stable category $\underline{\tt GProj}(kH)$ admits a set
${\tt S}_H$ of compact generators that consists of Gorenstein
projective $kH$-modules of type FP$_{\infty}$. Applying Corollary
4.3 to the case where the class $\mathcal{B}_0$ therein is the
subclass ${\mathfrak{G}_{AC}} \subseteq \mathfrak{B}$ (cf.\
Corollary 4.4), we conclude that the stable category
$\underline{\tt GProj}(kG)$ admits
\[ {\tt T}_G = {\textstyle{\bigcup}} \left\{
   \mbox{ind}_H^G {\tt S}_H : H \subseteq G ,
   H \in {\mathfrak{G}_{AC}} \right\} \]
as a set of compact generators. Since induction from subgroups
to $G$ preserves the modules of type FP$_{\infty}$, the generating
set ${\tt T}_G$ consists of $kG$-modules of type FP$_{\infty}$.
We now consider the triangulated subcategory
$\langle {\tt T}_G \rangle$ of the stable category
$\underline{\tt GProj}(kG)$ that is generated by ${\tt T}_G$.
Then, the full subcategory of compact objects of
$\underline{\tt GProj}(kG)$ is the thick closure
$\mbox{thick} \langle {\tt T}_G \rangle$ of
$\langle {\tt T}_G \rangle$ in $\underline{\tt GProj}(kG)$;
see, for example, \cite[Proposition 3.4.15]{Kra2}. It follows 
that $M \in \mbox{thick} \langle {\tt T}_G \rangle$. We now 
invoke Lemma 3.2 and conclude that
$\langle {\tt T}_G \rangle$ is contained in
$\underline{\tt GProj}^{fin}(kG)$. Using the notation 
following Lemma 3.2, we therefore obtain an inclusion
\[ \mbox{thick} \langle {\tt T}_G \rangle \subseteq
   \mbox{thick} \! \left[ \underline{\tt GProj}^{fin}(kG)
   \right] \! = {\tt K}(R) . \]
It follows that $M \in {\tt K}(R)$, so that Lemma 3.3 implies 
that $M$ is actually of type FP$_{\infty}$. \hfill $\Box$

\begin{Corollary}
Assume that $k$ is coherent and Gorenstein regular and consider
a group $G$ contained in $\overline{\mathfrak{F}}$. Then, the
following conditions are equivalent for a $kG$-module $M$:

(i) $M$ is a cokernel of an acyclic complex of finitely
    generated free $kG$-modules,

(ii) $M$ is Gorenstein projective of type FP$_{\infty}$,

(iii) $M$ is Gorenstein projective and finitely presented and

(iv) $M$ is Gorenstein projective and finitely generated.
\end{Corollary}
\vspace{-0.05in}
\noindent
{\em Proof.}
As shown in the proof of Corollary 4.5, the assumptions on
$k$ imply that $\mathfrak{G}_{AC}$ contains all finite groups,
i.e.\ that $\mathfrak{F} \subseteq \mathfrak{G}_{AC}$. Then,
we also have
$\overline{\mathfrak{F}} \subseteq \overline{\mathfrak{G}_{AC}}$
and the result follows from Proposition 5.1. \hfill $\Box$

\vspace{0.1in}

\noindent
{\bf Remark 5.3.}
The results presented in this Section are  proved using the
general theory of triangulated categories, even though their
statements make no reference whatsoever to stable or triangulated
categories. We are not aware of any proof of these results
that avoids the use of triangulated categories.

\section{The Gorenstein projective model structure in $kG$-Mod}

\noindent
In this section, we shall use the results obtained earlier
in the paper and demonstrate that the Gorenstein projective
model structure has a compactly generated tensor triangulated
homotopy category over many $\mathfrak{X}$-groups. In particular,
we shall prove the Theorem stated in the Introduction.

We consider a commutative ring $k$ and note that for any
$\mathfrak{X}$-group $G$ the Gorenstein projective cotorsion
pair $\left( {\tt GProj}(kG),{\tt GProj}(kG)^{\perp} \right)$
is complete, projective and hereditary; this is property (3)
in $\S $1.IV. Since ${\tt GProj}(kG)^{\perp}$ is a thick class,
we may consider the Hovey triple
\[ \left( {\tt GProj}(kG) , {\tt GProj}(kG)^{\perp} ,
   kG\mbox{-Mod} \right) \]
and the associated exact model structure $\mathcal{M}$ in the
category of $kG$-modules, defined in \cite{Ho2}. We can make
explicit this structure as follows:
\begin{itemize}
\item All $kG$-modules are fibrant in $\mathcal{M}$, whereas the
      cofibrant objects are the Gorenstein projective $kG$-modules
      and the trivial objects are the modules in the orthogonal
      class ${\tt GProj}(kG)^{\perp}$.
\item The fibrations (resp.\ the trivial fibrations) are the
      surjective $kG$-linear maps (resp.\ the surjective
      $kG$-linear maps with kernel in ${\tt GProj}(kG)^{\perp}$).
\item The cofibrations (resp.\ the trivial cofibrations) are the
      injective $kG$-linear maps whose cokernel is Gorenstein
      projective (resp.\ the injective $kG$-linear maps whose
      cokernel is projective).
\item As usual, the weak equivalences are those $kG$-linear maps
      that may be expressed as compositions of a trivial cofibration
      followed by a trivial fibration.
\end{itemize}
The homotopy category $\mbox{Ho}(\mathcal{M})$ of the model structure
$\mathcal{M}$, which is obtained from the module category by formally
inverting the weak equivalences, is a triangulated category which is
triangulated equivalent to the stable category $\underline{\tt GProj}(kG)$
of Gorenstein projective $kG$-modules. This follows from Gillespie's
fundamental result \cite[Proposition 4.4 and Corollary 4.8]{Gil1}.

Let $\mathfrak{A}$ and $\mathfrak{B}$ be the classes defined in
Sections 2 and 4 respectively and consider the class $\mathfrak{C}$
consisting of those groups that are contained in both $\mathfrak{A}$
and $\mathfrak{B}$. In view of Theorems 2.3 and 4.2, the class
$\mathfrak{C} = \mathfrak{A} \cap \mathfrak{B}$ is closed under the
operations ${\scriptstyle{{\bf LH}}}$, $\Phi$ and $\Phi_{flat}$.

\begin{Theorem}
Assume that $k$ has finite weak global dimension and let $G$
be a $\mathfrak{C}$-group. Then, the homotopy category
$\mbox{Ho}(\mathcal{M})$ of the Gorenstein projective model
structure $\mathcal{M}$ in the category of $kG$-modules is a
compactly generated tensor triangulated category.
\end{Theorem}
\vspace{-0.05in}
\noindent
{\em Proof.}
The tensor product of $kG$-modules (with the diagonal action of
$G$) endows the category $kG$-Mod with a closed symmetric monoidal
structure. The closed structure is provided by letting $G$ act
diagonally on the Hom-group $\mbox{Hom}_k(\_\!\_,\_\!\_)$ between
two $kG$-modules. We shall use \cite[Theorem 4.3.2]{Ho1}, in order
to show that the derived tensor product endows the homotopy category
$\mbox{Ho}(\mathcal{M})$ with the structure of a closed symmetric
model category. To that end, it suffices to show that $\mathcal{M}$
is a {\em monoidal} model structure in the category of $kG$-modules,
in the sense of \cite[Definition 4.2.6]{Ho1}. Invoking \cite[$\S $7]{Ho2},
it therefore suffices to verify that the following conditions are
satisfied:

(i) All cofibrations are $k$-pure monomorphisms.

(ii) The tensor product of Gorenstein projective $kG$-modules
is Gorenstein projective.

(iii) The class of projective $kG$-modules is invariant under
tensoring with Gorenstein projective $kG$-modules.

(iv) If $N$ is a Gorenstein projective $kG$-module, then a special
Gorenstein projective precover $f : M \longrightarrow k$ induces a
weak equivalence
$1 \otimes f : N \otimes_kM \longrightarrow k \otimes_kN=N$.
\newline
Indeed, (iv) is precisely condition 2 in \cite[Definition 4.2.6]{Ho1},
whereas (i), (ii) and (iii) imply condition 1 therein; see the proof
of \cite[Theorem 7.2]{Ho2}.

Our assumption on $k$ implies that all Gorenstein projective $kG$-modules
are $k$-projective; see the proof of Proposition 2.1. In particular, all
Gorenstein projective $kG$-modules are $k$-flat and hence (i) follows.
Assertion (ii) is precisely Proposition 2.1, whereas (iii) follows since
the tensor product of a projective $kG$-module with a $k$-projective
$kG$-module is also a projective $kG$-module. As far as (iv) is concerned,
we note that $\mbox{ker} \, f \in {\tt GProj}(kG)^{\perp}$. Since
$G \in \mathfrak{C} \subseteq \mathfrak{A}$ and $N \in {\tt GProj}(kG)$,
we conclude that
$\mbox{ker} \, (1 \otimes f) = N \otimes_k \mbox{ker} \, f \in
 {\tt GProj}(kG)^{\perp}$.
It follows that $1 \otimes f$ is a trivial fibration and hence a weak
equivalence; this proves assertion (iv).

Finally, we note that the homotopy category of $\mathcal{M}$ is
compactly generated, since
$\mbox{Ho}(\mathcal{M}) \simeq \underline{\tt GProj}(kG)$ and
$G \in \mathfrak{C} \subseteq \mathfrak{B}$. \hfill $\Box$

\vspace{0.1in}

\noindent
The following result extends \cite[Theorem 1.1]{Ke}.

\begin{Corollary}
Assume that $k$ has finite weak global dimension and let $G$
be a $\overline{\mathfrak{G}_{AC}}$-group. Then, the homotopy
category $\mbox{Ho}(\mathcal{M})$ of the Gorenstein projective
model structure $\mathcal{M}$ in the category of $kG$-modules
is a compactly generated tensor triangulated category.
\end{Corollary}
\vspace{-0.05in}
\noindent
{\em Proof.}
We note that $\mathfrak{G}_{AC} \subseteq \mathfrak{G}$,
where $\mathfrak{G}$ is the class consisting of those
groups $G$, for which $kG$ is Gorenstein regular. In
particular, $\mathfrak{G}_{AC} \subseteq \mathfrak{W}$,
where $\mathfrak{W}$ is the class consisting of those
groups $G$, for which $kG$ is weakly Gorenstein regular.
Since $k$ has finite weak global dimension, Corollary 2.4
implies that
$\overline{\mathfrak{G}_{AC}} \subseteq \mathfrak{A}$.
On the other hand, Corollary 4.4 implies that
$\overline{\mathfrak{G}_{AC}} \subseteq \mathfrak{B}$, so
that
$\overline{\mathfrak{G}_{AC}} \subseteq
 \mathfrak{A} \cap \mathfrak{B} = \mathfrak{C}$.
Hence, the result follows from Theorem 6.1.
\hfill $\Box$

\vspace{0.1in}

\noindent
We now consider the closure $\overline{\mathfrak{F}}$ of the
class $\mathfrak{F}$ of finite groups under the operation
${\scriptstyle{{\bf LH}}} \cup \Phi \cup \Phi_{flat}$ and
prove the following result, which is precisely the Theorem
stated in the Introduction.

\begin{Corollary}
Assume that $k$ is a coherent ring of finite global dimension
and let $G$ be an $\overline{\mathfrak{F}}$-group. Then, the
homotopy category $\mbox{Ho}(\mathcal{M})$ of the Gorenstein
projective model structure $\mathcal{M}$ in the category of
$kG$-modules is a compactly generated tensor triangulated
category.
\end{Corollary}
\vspace{-0.05in}
\noindent
{\em Proof.}
Since $k$ has finite weak global dimension, Corollary 2.5
implies that $\overline{\mathfrak{F}} \subseteq \mathfrak{A}$.
On the other hand, $k$ is coherent and Gorenstein
regular\footnote{A ring has finite global dimension if and
only if it is Gorenstein regular and has finite weak global
dimension.} and hence
$\overline{\mathfrak{F}} \subseteq \mathfrak{B}$; cf.
Corollary 4.5. We conclude that
$\overline{\mathfrak{F}} \subseteq
 \mathfrak{A} \cap \mathfrak{B} = \mathfrak{C}$,
so that the result is again a consequence of Theorem 6.1.
\hfill $\Box$

\appendix

\section{Groups outside $\overline{\mathfrak{F}}$}

\noindent
Our goal in this appendix is to show that there are many finitely
generated groups that are not contained in the closure of the class
of finite groups under the operation
${\scriptstyle{{\bf LH}}} \cup \Phi \cup \Phi_{flat}$.

Let $k$ be a commutative ring. It has been shown by Fournier-Facio
and Sun in \cite[Theorem 5.7]{FS} that there exist uncountably many
mutually non-isomorphic finitely generated torsion-free simple groups
of infinite homological and cohomological dimension over $k$, all of
whose non-trivial proper subgroups are infinite cyclic. Any cellular
action of such a group on a finite dimensional contractible CW-complex
has necessarily a fixed point; cf.\ \cite[Corollary 5.8]{FS}. We fix
such a group $G$ below and let $\mathfrak{D}$ be any class of groups.

\medskip

\noindent
(i) If $G$ is an ${\scriptstyle{{\bf LH}}}\mathfrak{D}$-group, then
$G \in \mathfrak{D}$.
\newline
{\em Proof.}
Assume that $G \in {\scriptstyle{{\bf LH}}}\mathfrak{D}$. Since $G$
is finitely generated, it is an ${\scriptstyle{{\bf H}}}\mathfrak{D}$-group.
Let $\alpha$ be the smallest ordinal with
$G \in {\scriptstyle{{\bf H}}}_{\alpha}\mathfrak{D}$. We claim that
$\alpha =0$, so that
$G \in {\scriptstyle{{\bf H}}}_0\mathfrak{D} = \mathfrak{D}$,
as needed. Assume on the contrary that $\alpha >0$. Then, $\alpha$
is a successor ordinal, say $\alpha = \beta +1$. Since $G$ is an
${\scriptstyle{{\bf H}}}_{\beta +1}\mathfrak{D}$-group, it acts
cellularly on a finite dimensional contractible CW-complex with
all isotropy subgroups contained in
${\scriptstyle{{\bf H}}}_{\beta}\mathfrak{D}$. This action must
fix a cell, so that $G$ itself is such an isotropy subgroup. It
follows that $G \in {\scriptstyle{{\bf H}}}_{\beta}\mathfrak{D}$,
contradicting the minimality of $\alpha$. \hfill $\Box$

\medskip

\noindent
(ii) If $G$ is a $\Phi \mathfrak{D}$-group, then $G \in \mathfrak{D}$.
\newline
{\em Proof.}
Assume on the contrary that $G \in \Phi \mathfrak{D}$ but
$G \notin \mathfrak{D}$ and consider the trivial $kG$-module $k$.
If $H$ is a $\mathfrak{D}$-subgroup of $G$, then $H \neq G$ (since
$G \notin \mathfrak{D}$) and hence $H$ is either trivial or infinite
cyclic; in any case, $\mbox{pd}_{kH}k \leq 1$. Since
$G \in \Phi \mathfrak{D}$, it follows that $\mbox{pd}_{kG}k < \infty$.
This is a contradiction, since $G$ has infinite cohomological dimension
over $k$. \hfill $\Box$

\medskip

\noindent
(iii) If $G$ is a $\Phi_{flat} \mathfrak{D}$-group, then $G \in \mathfrak{D}$.
\newline
{\em Proof.}
The argument is completely analogous to that given in (ii) above.
\hfill $\Box$

\medskip

\noindent
Recall from $\S $1.II that we may define for any ordinal $\alpha$
the class $\mathfrak{D}_{\alpha}$ by transfinite induction, letting
$\mathfrak{D}_0 = \mathfrak{D}$,
$\mathfrak{D}_{\alpha} =
 {\scriptstyle{{\bf LH}}}\mathfrak{D}_{\beta} \cup
 \Phi\mathfrak{D}_{\beta} \cup \Phi_{flat}\mathfrak{D}_{\beta}$
if $\alpha = \beta +1$ and
$\mathfrak{D}_{\alpha} = \bigcup_{\beta < \alpha}\mathfrak{D}_{\beta}$
if $\alpha$ is a limit ordinal. Then, the closure $\overline{\mathfrak{D}}$
of $\mathfrak{D}$ under the operation
${\scriptstyle{{\bf LH}}} \cup \Phi \cup \Phi_{flat}$ is the class
consisting of those groups which are contained in $\mathfrak{D}_{\alpha}$,
for some ordinal $\alpha$.

\medskip

\noindent
(iv) If $G$ is a $\overline{\mathfrak{D}}$-group, then
$G \in \mathfrak{D}$.
\newline
{\em Proof.}
Let $\alpha$ be the smallest ordinal with
$G \in \mathfrak{D}_{\alpha}$. We claim that $\alpha =0$, so that
$G \in \mathfrak{D}_0 = \mathfrak{D}$. Assume on the contrary that
$\alpha >0$. Then, $\alpha$ is a successor ordinal, say
$\alpha = \beta +1$. Since $G$ is contained in
$\mathfrak{D}_{\alpha} =
 {\scriptstyle{{\bf LH}}}\mathfrak{D}_{\beta} \cup
 \Phi\mathfrak{D}_{\beta} \cup \Phi_{flat}\mathfrak{D}_{\beta}$,
we may use assertions (i), (ii) and (iii) above (applied to the class
$\mathfrak{D}_{\beta}$) and conclude that $G \in \mathfrak{D}_{\beta}$.
This is a contradiction, in view of the minimality of $\alpha$.
\hfill $\Box$

\begin{Proposition}
The groups constructed in \cite[Theorem 5.7]{FS} are not contained in
$\overline{\mathfrak{F}}$.
\end{Proposition}
\vspace{-0.05in}
\noindent
{\em Proof.}
This is an immediate consequence of assertion (iv) above (by
letting $\mathfrak{D}$ be the class $\mathfrak{F}$ of finite
groups therein), since the groups constructed in
\cite[Theorem 5.7]{FS} are infinite. Indeed, these groups are
torsion-free and have infinite (co-)homological dimension.
\hfill $\Box$

\vspace{0.1in}

\noindent
{\em Acknowledgments.}
It is a pleasure to thank Apostolos Beligiannis, Gregory
Kendall and Wei Ren for useful discussions on the topic
of this paper.

\vspace{0.1in}

\noindent
{\small {\sc Department of Mathematics,
             University of Athens,
             Athens 15784,
             Greece}}

\noindent
{\em E-mail addresses:} {\tt emmanoui@math.uoa.gr} and
                        {\tt otalelli@math.uoa.gr}

\end{document}